\RequirePackage{fix-cm}
\documentclass[smallextended]{svjour3}
\smartqed  % flush right qed marks, e.g. at end of proof
\usepackage{graphics}
\usepackage{epsfig}
\usepackage{amssymb}
\usepackage{amsfonts}
\usepackage{amsmath}
\usepackage{graphicx}
\journalname{}
\begin{document}
\title{Stepanov-like Weighted Pseudo-Almost Automorphic Solutions on Time Scales for a Novel High-order BAM Neural Network with Mixed
Time-varying Delays in the Leakage Terms}
\author{Adn\`ene Arbi}
\institute{
Adn\`ene Arbi\\
adnen.arbi@enseignant.edunet.tn, arbiadnene@yahoo.fr\\
Laboratory of Engineering Mathematics (LR01ES13), Tunisia
Polytechnic School, University of Carthage, El Khawarizmi Street, Carthage 2078, Tunisia.\\
} \maketitle
\begin{abstract}
We first propose the concept of Stepanov-like weighted pseudo almost automorphic on time-space scales and we apply this type of oscillation
to high-order BAM neural networks with mixed delays. Then, we study the existence and exponential stability of Sp-weighted pseudo-almost automorphic
on time-space scales solutions for the suggested system. Some criteria assuring the convergence are proposed. Our method is mainly based on the Banachs
fixed point theorem, the theory of calculus on time scales and the Lyapunov- Krasovskii functional method. Moreover, a numerical example is given to show
the effectiveness of the main results.
\end{abstract}
\vspace{-6pt}
\section{Introduction}
The concept of weighted pseudo almost automorphic on time-space
scales functions for the nabla and delta derivative is recently
introduced (see \cite{tomorphic2}, \cite{chang},
\cite{adnene+ahmed}). It is a natural generalization of almost
automorphic on time-space scales functions introduced in
\cite{carlos}. In 2010, the concept of Stepanov-like weighted pseudo
almost automorphy which is a natural generalization of the almost
automorphy is presented (see \cite{diagana}). Moreover, there is no
definition of the notion of Stepanov-like almost automorphy and
Stepanov-like weighted pseudo almost periodic on time-space scales
in the previous work.

On the other hand, many researchers have been devoting the dynamics
of various class of neural networks due to its wide application in
pattern recognition, associative memory, image, and signal
processing (see \cite{adnene+ahmed}, \cite{add1}, \cite{add2}, \cite{adn0}, \cite{adn00}, \cite{adn01}, \cite{adn1},
\cite{aaa}, \cite{adnene+jinde}, \cite{maas}, \cite{kren},
\cite{zhang1}, \cite{zhang2}, \cite{zhang3}). Furthermore, BAM
neural networks as an extension of the unidirectional autoassociator
of Hopfield neural network (see \cite{adn0},
\cite{adn01}, \cite{adn1}, \cite{aaa}), was firstly introduced by
Kosko \cite{kosko}. In recent years, many scholars pay much
attention to the dynamical behavior of bidirectional associative
memory (BAM) neural networks. Considering that time delays are
unavoidable because of the finite switching of amplifiers in
practical implementation of neural networks, and the time delay may
result in oscillation and instability; many authors focus on the
dynamical properties of BAM neural networks with time delays (see
\cite{maas}, \cite{bam0}, \cite{bam1}, \cite{bam2}, \cite{bam3},
\cite{bam4}, \cite{bam5}, \cite{bam6}). In real application, when
robot move, the joints are properly described by almost periodic
solutions of a dynamic neural network. For this reason, it is very
important to study almost periodic solutions of neural networks
models. In \cite{bam0}, the authors investigated the almost periodic
solution of
\begin{equation*} \left\{
\begin{array}{ll}
\dot{x}_{i}(t)=-\alpha_{i}(t)x_{i}(t)+\sum\limits_{j=1}^{m}b_{ij}(t)(t)f_{j}(y_{j}(t-\tau))+I_{i}(t), \\
\dot{y}_{j}(t)=-c_{j}(t)y_{j}(t)+\sum\limits_{j=1}^{n}\bar{b}_{ij}(t)(t)f_{j}(x_{j}(t-\sigma))+J_{j}(t),\,\
t\geq 0\\
x_{i}(t)=\phi_{i}(t), \,\ y_{j}(t)=\psi_{j}(t), \,\
t\in[-\tau^{*},0], \\ i=1,...,n \,\ j=1,...,m,  \ \ n,m\in
\mathbb{Z}_{+}. \text{}
\end{array}\right.
\end{equation*}

The generalized high-order BAM neural network with mixed delays,
defined as follows, has faster convergence rate, higher fault
tolerance, and stronger approximation property. The problem of
existence and global exponential stability of periodic solution for
high-order discrete-time BAM neural networks has been studied in
\cite{highbam1}. In fact, the study of the existence periodic
solutions, as well as its numerous generalizations to almost
periodic solutions, pseudo almost periodic solutions, weighted
pseudo almost periodic solutions, and so forth, is one of the most
attracting topics in the qualitative theory of differential
equations due both to its mathematical interest as well as to their
applications in various areas of applied science. The authors in
\cite{highbam} propose some several sufficient conditions for
ensuring existence, global attractivity and global asymptotic
stability of the periodic solution for the higher-order
bidirectional associative memory neural networks with periodic
coefficients and delays by using the continuation theorem of
Mawhin's coincidence degree theory, the Lyapunov functional and the
non-singular $M$-matrix:
\begin{equation*} \left\{
\begin{array}{ll}
\dot{x}_{i}(t)=-\alpha_{i}x_{i}(t)+\sum\limits_{j=1}^{m}b_{ij}(t)f_{j}(y_{j}(t-\tau))
\\+\sum\limits_{j=1}^{m}\sum\limits_{l=1}^{m}c_{ijl}(t)f_{j}(y_{j}(t-\tau))f_{l}(y_{l}(t-\tau))
+I_{i}(t),
\\
\dot{y}_{j}(t)=-c_{j}y_{j}(t)+\sum\limits_{j=1}^{n}\bar{b}_{ij}(t)f_{j}(x_{j}(t-\sigma))
\\+\sum\limits_{j=1}^{n}\sum\limits_{l=1}^{n}\bar{c}_{ijl}(t)g_{i}(x_{i}(t-\sigma))g_{l}(x_{l}(t-\sigma))
+J_{j}(t),\,\ t\geq 0.
\end{array}\right.
\end{equation*}
With initial condition
\begin{equation*} \left\{
\begin{array}{ll}
x_{i}(t)=\phi_{i}(t), \,\ y_{j}(t)=\psi_{j}(t), \,\
t\in[-\tau^{*},0],\\ i=1,...,n \,\ j=1,...,m,  \ \ n,m\in
\mathbb{Z}_{+}. \text{}
\end{array}\right.
\end{equation*}

Furthermore, it has been reported that if the parameters and time
delays are appropriately chosen, the delayed high-order BAM neural
network can exhibit complicated behaviors even with strange chaotic
attractors. Based on the aforementioned arguments, the study of the
high-order BAM neural network with mixed delays and its analogous
equations have attracted worldwide interest (see \cite{bam6},
\cite{ada1}, \cite{ada2}, \cite{ada3}, \cite{ada4}, \cite{highbam1},
\cite{highbam}, \cite{leak1}, \cite{leak2}, \cite{leak3}). In fact,
it is important that systems contain some information about the
derivative of the past state to further describe the dynamics for
such complex neural reactions. In real world, the mixed time-varying
delays and leakage delay should be taken into account when modeling
realistic neural networks (see \cite{adnene+ahmed},
\cite{adnene+jinde}, \cite{maas}).

As a continuation of our previous published results, we shall
consider a high-order BAM neural network with mixed delays:
\begin{equation*}
\left\{
\begin{array}{ccc}
x^{\Delta}_{i}\left( t\right) =
-\alpha_{i}(t)x_{i}\left(t-\eta_{i}(t)\right)
+\sum\limits_{j=1}^{m}D_{ij}\left( t\right) f_{j}\left(x_{j}\left(
t\right) \right) \\+\sum\limits_{j=1}^{m}D_{ij}^{\tau}\left(
t\right) f_{j}\left(x_{j}\left( t -\tau_{ij}(t)\right) \right)
+\sum\limits_{j=1}^{m}\overline{D}_{ij}\left(
t\right)\int_{t-\sigma_{ij}(t)} ^{t}f_{j}\left(x_{j}\left( s\right)
\right)\Delta s \\+\sum\limits_{j=1}^{m}\widetilde{D}_{ij}\left(
t\right)\int_{t-\xi_{ij}(t)} ^{t}f_{j}\left(x_{j}^{\Delta}\left(
s\right) \right)\Delta s+I_{i}\left( t\right)\\
+\sum\limits_{j=1}^{m}\sum\limits_{k=1}^{m}T_{ijk}(t)f_{k}(x_{k}(t-\chi_{k}(t)))f_{j}(x_{j}(t-\chi_{j}(t))),\,\ i=1,...,n, \\
y^{\Delta}_{j}\left( t\right) =
-c_{j}(t)y_{j}\left(t-\varsigma_{j}(t)\right)
+\sum\limits_{i=1}^{n}E_{ij}\left( t\right) f_{j}\left(x_{j}\left(
t\right) \right) \\+\sum\limits_{i=1}^{n}E_{ij}^{\tau}\left(
t\right) f_{j}\left(x_{j}\left( t -\tau_{ij}(t)\right) \right)
+\sum\limits_{i=1}^{n}\overline{E}_{ij}\left(
t\right)\int_{t-\sigma_{ij}(t)} ^{t}f_{j}\left(x_{j}\left( s\right)
\right)\Delta s \\+\sum\limits_{i=1}^{n}\widetilde{E}_{ij}\left(
t\right)\int_{t-\xi_{ij}(t)}^{t}f_{j}\left(x_{j}^{\Delta}\left(
s\right) \right)\Delta s+J_{j}\left( t\right)\\
+\sum\limits_{i=1}^{n}\sum\limits_{k=1}^{n}\overline{T}_{ijk}(t)f_{k}(x_{k}(t-\chi_{k}(t)))f_{j}(x_{j}(t-\chi_{j}(t)))
,\,\ t\in\mathbb{T}, \,\ j=1,...,m.
\end{array}
\right.
\end{equation*}

To the best of our knowledge, the existence of Stepanov-like
weighted pseudo-almost automorphic solution on time-space scales to
BAM neural networks (BAMs) and high-order BAM neural network
(HOBAMs) with variable coefficients, mixed delays and leakage delays
have not been studied yet. It has been reported that if the
parameters and time delays are appropriately chosen, the delayed
neural networks in leakage term can exhibit complicated behaviors
even with strange chaotic attractors (see \cite{adnene+ahmed},
\cite{adnene+jinde}, \cite{maas}, \cite{leak1}, \cite{leak2},
\cite{leak3}, \cite{leak0}). In addition, the theory of time scales,
which has recently received much attention, was introduced by Hilger
in his PhD thesis in 1988 to unify continuous and discrete analysis
\cite{hilger}.

Our main purpose of this paper is to introduce the Stepanov-like
weighted pseudo almost automorphic functions on time-space scales,
study some of their basic properties and establish the existence,
uniqueness, stability and convergence of Stepanov-like weighted
pseudo almost automorphic solutions of HOBAMs on time-space scales.
we prove new composition theorems for Stepanov-like weighted pseudo
almost automorphic functions on time-space scales.

The remainder of this paper is organized as follows. In Section 2,
we will present the model of HOBAMs. In section 3, we will introduce
some necessary notations, definitions and fundamental properties of
the weighted pseudo-almost automorphic on time-space scales
environment, which will be used in the paper. In Section 4, some
sufficient conditions will be derived ensuring the existence of the
Stepanov-like weighted pseudo-almost automorphic solution on
time-space scales. Section 5 will be devoted to the exponential
stability of the Stepanov-like weighted pseudo-almost automorphic
solution on time-space scales of a HOBAMs model, and the convergence
of all solutions to its unique Stepanov-like weighted pseudo-almost
automorphic solution. At last, one illustrative numerical example
will be given.
\section{Preliminaries and function spaces}
In the following, we introduce some definitions and state some
preliminary results.
\subsection{Time-space scales and delta derivative}
\begin{definition}(\cite{adnene+ahmed})
Let $\mathbb{T}$ be a nonempty closed subset (time scale) of
$\mathbb{R}$. The forward and backward jump operators $\sigma,\rho:
\mathbb{T}\longrightarrow \mathbb{T}$ and the graininess
$\nu:\mathbb{T}\longrightarrow\mathbb{R}_{+}$ are defined,
respectively, by $\sigma(t)=\inf\{s\in\mathbb{T}:s>t\}$,
$\rho(t)=\sup\{s\in\mathbb{T}:s<t\}$ and $\nu(t)=\sigma(t)-t$.
\end{definition}
\begin{lemma}(\cite{adnene+ahmed}, \cite{advance})
Considering that $f,g$ be delta differentiable functions on
$\mathbb{T}$, then:
\begin{description}
\item[(i)] $(\lambda_{1}f+\lambda_{2}g)^{\Delta}=\lambda_{1}f^{\Delta}+\lambda_{2}g^{\Delta}$, for any constants $\lambda_{1},\lambda_{2}$;
\item[(ii)] $(fg)^{\Delta}(t)=f^{\Delta}(t)g(t)+f(\sigma(t))g^{\Delta}(t)=f(t)g^{\Delta}(t)+f^{\Delta}(t)g(\sigma(t))$;
\item[(iii)] If $f$ and $f^{\Delta}$ are continuous, then $\left(\int_{a}^{t}f(t,s)\Delta s\right)^{\Delta}=f(\sigma(t),t)+\int_{a}^{t}f(t,s)\Delta s$.
\end{description}
\end{lemma}
\begin{lemma}(\cite{adnene+ahmed}, \cite{advance})  Assume that $p,q:\mathbb{T}\longrightarrow\mathbb{R}$ are two regressive functions, then
\begin{description}
\item[(i)] $e_{0}(t,s)\equiv1$ and $e_{p}(t,t)\equiv1$;
\item[(ii)] $e_{p}(t,s)=\frac{1}{e_{p}(s,t)}=e_{\ominus p}(s,t)$;
\item[(iii)] $e_{p}(t,s)e_{p}(s,r)=e_{p}(t,r)$;
\item[(iv)] $\left(e_{p}(t,s)\right)^{\Delta}=p(t)e_{p}(t,s)$.
\end{description}
\end{lemma}
\begin{lemma}(\cite{adnene+ahmed}, \cite{advance})
Assume that $p(t)\geq0$ for $t\geq s$, then $e_{p}(t,s)\geq1$.
\end{lemma}
\begin{definition}(\cite{adnene+ahmed})
A function $p:\mathbb{T}\longrightarrow\mathbb{R}$ is called
regressive provided $1+\mu(t)p(t)\neq0$ for all
$t\in\mathbb{T}^{k}$; $p:\mathbb{T}\longrightarrow\mathbb{R}$ is
called positively provided $1+\mu(t)p(t)>0$ for all
$t\in\mathbb{T}^{k}$. The set of all regressive and rd-continuous
functions $p:\mathbb{T}\longrightarrow\mathbb{R}$ will be denoted by
$\mathcal{R}=\mathcal{R}(\mathbb{T},\mathbb{R})$ and the set of all
positively regressive functions and rd-continuous functions will be
denoted $\mathcal{R}^{+}=\mathcal{R}^{+}(\mathbb{T},\mathbb{R})$.
\end{definition}
\begin{lemma} (\cite{adnene+ahmed}, \cite{advance})
Suppose that $p\in\mathcal{R}^{+}$, then:
\begin{description}
\item[(i)] $e_{p}(t,s)>0$, for all $t,s\in\mathbb{T}$;
\item[(ii)] if $p(t)\leq q(t)$ for all $t\geq s$, $t,s\in\mathbb{T}$, then $e_{p}(t,s)\leq e_{q}(t,s)$ for all $t\geq s$.
\end{description}
\end{lemma}
\begin{lemma} (\cite{adnene+ahmed}, \cite{advance})
If $p\in\mathcal{R}$ and $a,b,c\in \mathbb{T}$, then $[e_{p}(c,.)]^{\Delta}=-p[e_{p}(c,.)]^{\sigma}$,\\
and $\int_{a}^{b}p(t)e_{p}(c,\sigma(t))\Delta
t=e_{p}(c,a)-e_{p}(c,b)$.
\end{lemma}
\begin{lemma} (\cite{adnene+ahmed}, \cite{advance})
Let $a\in\mathbb{T}^{k}$, $b\in\mathbb{T}$ and assume that
$f:\mathbb{T}\times\mathbb{T}^{k}\longrightarrow\mathbb{R}$ is
continuous at $(t,t)$, where $t\in\mathbb{T}^{k}$ with $t>a$.
Additionally assume that $f^{\Delta}(t,.)$ is  rd-continuous on
$[a,\sigma(t)]$. Suppose that for each $\epsilon>0$, there exists a
neighborhood $U$ of $\tau\in[a,\sigma(t)]$ such that
\begin{equation*}
|f(\sigma(t),\tau)-f(s,\tau)-f^{\Delta}(t,\tau)(\sigma(t)-s)|\leq\epsilon|\sigma(t)-s|,
\,\ \forall s\in U,
\end{equation*}
where $f^{\Delta}$ denotes the derivative of $f$ with respect to the
first variable. Then
\begin{description}
\item[(i)] $g(t):=\int_{0}^{t}f(t,\tau)\Delta\tau$ implies $g^{\Delta}(t):=\int_{a}^{t}f^{\Delta}(t,\tau)\Delta\tau+f(\sigma(t),t)$;
\item[(ii)] $h(t):=\int_{t}^{b}f(t,\tau)\Delta\tau$ implies $h^{\Delta}(t):=\int_{t}^{b}f^{\Delta}(t,\tau)\Delta\tau-f(\sigma(t),t)$.
\end{description}
\end{lemma}

For more details of time scales and $\Delta$-measurability, one is
referred to read the excellent books (\cite{advance},
\cite{guseinov}).
\subsection{Stepanov-like weighted pseudo-almost automorphic functions on time-space scales}
In the following, we recall some definitions of Stepanov almost
automorphic functions and Stepanov-like weighted pseudo-almost
automorphic functions on time-space scales.
\begin{definition}(\cite{tomorphic2})
A time scale $\mathbb{T}$ is called an almost periodic time scale if
\begin{equation*}
\Pi:=\left\{\tau\in\mathbb{R}: t\pm\tau\in\mathbb{T}, \forall
t\in\mathbb{T}\right\}\neq0.
\end{equation*}
\end{definition}
\begin{definition}
A function $f:\mathbb{T}\longrightarrow\mathbb{R}$ is Bochner
integrable, or integrable for short, if there is a sequence of
functions such that $f_{n}(t)\longrightarrow f(t)$ pointwise a.e. in
$\mathbb{T}$ and
\begin{equation*}
\lim\limits_{n\longrightarrow+\infty}\int_{\mathbb{T}}\|f(s)-f_{n}(s)\|\Delta
s=0,
\end{equation*}
and the integral of $f$ is defined by
\begin{equation*}
\int_{\mathbb{T}}f(s)\Delta
s=\lim\limits_{n\longrightarrow+\infty}f_{n}(s)\Delta s,
\end{equation*}
where the limit exists strongly in $\mathbb{R}$.
\end{definition}
\begin{definition}
Let $E\subset\mathbb{T}$ and $f:\mathbb{T}\longrightarrow\mathbb{R}$
be a strongly $\Delta$-measurable function. If, for a given $p$,
$1\leq p<\infty$, $f$ satisfies
\begin{equation*}
\int_{K}\|f(s)\|^{p}\Delta s<+\infty,
\end{equation*}
where $K$ is a compact subset of $E$, then $f$ is called
$p$-integrable in the Bochner sense. The set of all such functions
is denoted by $L_{loc}^{p}(\mathbb{T},\mathbb{R})$.
\end{definition}

From now on, for $a,b\in\mathbb{R}$ and $a\leq b$, we denote
$a^{*}=\inf\{s\in\mathbb{T}, s\geq a\}$,
$b^{*}=\inf\{s\in\mathbb{T}, s\geq b\}$ and for integrable function
$f$, we denote $\int_{a}^{b}f(s)\Delta
s=\int_{a^{*}}^{b^{*}}f(s)\Delta s$. Obviously,
$a^{*},b^{*}\in\mathbb{T}$. If $a\in\mathbb{T}$, then $a^{*}=a$; If
$b\in\mathbb{T}$, then $b^{*}=b$.

Throughout the rest of paper we fix $p$, $1\leq p<\infty$. We say
that a function $f\in L_{loc}^{p}(\mathbb{T},\mathbb{R})$ is
$p$-Stepanov bounded ($S_{l}^{p}$-bounded) if
\begin{equation*}
\|f\|_{S_{l}^{p}}=\sup\limits_{t\in\mathbb{T}}\left(\frac{1}{l}\int_{t}^{t+l}\|f(s)\|^{p}\Delta
s\right)^{\frac{1}{p}}<+\infty,
\end{equation*}
where $l>0$ is a constant.

We denote by $L_{S}^{p}$ the set of all $S_{l}^{p}$-bounded
functions from $\mathbb{T}$ into $\mathbb{R}$.

\begin{definition}
Let $\mathbb{T}$ be an almost periodic time scale. A function
$f(t):\mathbb{T}\longrightarrow \mathbb{R}^{n}$ is said to be
$S^{p}$-almost automorphic, if for any sequence
$\{s_{n}\}_{n=1}^{\infty}\subset\Pi$, there is a subsequence
$\{\tau_{n}\}_{n=1}^{\infty}\subset\{s_{n}\}_{n=1}^{\infty}$ such
that
\begin{equation*}
\|g(t)-f(t+\tau_{n})\|_{S_{l}^{p}}\longrightarrow0, \,\ \text{as}
\,\ n\longrightarrow+\infty,
\end{equation*}
is well defined for each $t\in\Pi$ and
\begin{equation*}
\|g(t-\tau_{n})-f(t)\|_{S_{l}^{p}}\longrightarrow0, \,\ \text{as}
\,\ n\longrightarrow+\infty,
\end{equation*}
for each $t\in\Pi$. Denote by $S^{p}AA(\mathbb{T},\mathbb{R}^{n})$
the set of all such functions.
\end{definition}

Let
\begin{equation*}
S^{p}AA(\mathbb{T})=\left\{f\in
S^{p}C_{rd}(\mathbb{T},\mathbb{R}^{n}): f \,\ \text{is Stepanov
almost automorphic}\right\}
\end{equation*}
and $S^{p}BC(\mathbb{T},\mathbb{R}^{n})$ denote the space of all
bounded continuous functions, in the Stepanov sens, from
$\mathbb{T}$ to $\mathbb{R}^{n}$.

Let $\mathcal{U}$ be the set of all functions
$\nu:\mathbb{T}\longrightarrow(0,+\infty)$ which are positive and
locally $\Delta$-integrable over $\mathbb{T}$. For given
$r\in(0,+\infty)\cap \Pi$, set \begin{equation}\label{poid}
m(r,\nu,t_{0}):=\int_{Q_{r}}\nu(s)\Delta s, \,\ \text{for each} \,\
\nu\in\mathcal{U},
\end{equation}
where $Q_{r}:=[t_{0}-r,t_{0}+r]_{\mathbb{T}}$
($t_{0}=\min\{[0,\infty)_{\mathbb{T}}\}$). If $\nu(s)=1$ for each
$\nu\in\mathbb{T}$, then
$\lim\limits_{t\longrightarrow\infty}\nu(Q_{r})=\infty$.
Consequently, we define the space of weights $\mathcal{U}_{\infty}$
by
\begin{equation*}
\mathcal{U}_{\infty}:=\left\{\nu\in\mathcal{U}: \,\
\lim\limits_{r\longrightarrow+\infty}m(r,\nu,t_{0})=+\infty\right\}.
\end{equation*}
In addition to the aforesaid section, we define the set of weights
$\mathcal{U}_{B}$ by
\begin{equation*}
\mathcal{U}_{B}:=\left\{\nu\in\mathcal{U}_{\infty}: \,\ \nu \,\
\text{is bounded in the Stepanov sens and }\,\
\inf\limits_{s\in\mathbb{T}}\nu(s)>0\right\}.
\end{equation*}
It is clear that $\mathcal{U}_{B}\subset
\mathcal{U}_{\infty}\subset\mathcal{U}$. Let
$S^{p}BCU^{(0)}(\mathbb{T},\mathbb{R}^{n})$ denote the space of all
bounded uniformly continuous functions from $\mathbb{T}$ to
$\mathbb{R}^{n}$,
\begin{multline*}
\phantom{+++++}S^{p}AA^{(0)}(\mathbb{T})=S^{p}AA^{(0)}(\mathbb{T},\mathbb{R}^{n})=\left\{f\in
S^{p}BCU(\mathbb{T},\mathbb{R}^{n}):\right.\\\left. f \,\ \text{is
Stepanov almost automorphic}\right\}\phantom{+++++}
\end{multline*}
and define for $t_{0}\in\mathbb{T}, r\in\Pi$, the class of functions
$WPAA_{0}(\mathbb{T},\nu,t_{0})$ as follows:
%It is clear that $\mathcal{U}_{B}\subset \mathcal{U}_{\infty}\subset\mathcal{U}$. Now, for $\nu\in\mathcal{U}_{\infty}$ define
\begin{multline*}
WPAA_{0}(\mathbb{T},\nu,t_{0})=\left\{f\in S^{p}BCU(\mathbb{T},\mathbb{R}^{n}): f \,\ \text{is delta measurable such that}\right.\\
\left.\lim\limits_{r\longrightarrow+\infty}\frac{1}{m(r,\nu,t_{0})}\int_{t_{0}-r}^{t_{0}+r}|f(s)|\nu(s)\Delta
s=0 \right\}.
\end{multline*}
We are now ready to introduce the sets $S^{p}WPAA(\mathbb{T},\nu)$
of Stepanov-like weighted pseudo-almost automorphic on time-space
scales functions:
\begin{definition}\label{defwpap}
A function $f\in S^{p}C_{rd}(\mathbb{T},\mathbb{R}^{n})$ is called
Stepanov-like weighted pseudo almost automorphic on time-space
scales if $f=g+\phi$, where $g\in S^{p}AA(\mathbb{T})$ and $\phi\in
WPAP_{0}(\mathbb{T},\nu)$. Denote by $S^{p}WPAA(\mathbb{T})$, the
set of Stepanov-like weighted pseudo-almost automorphic on
time-space scales functions.
\end{definition}
\subsection{Results on composition theorems}
By Definition \ref{defwpap}, one can easily show that
\begin{lemma}
Let $\phi\in S^{p}BC_{rd}(\mathbb{T},\mathbb{R}^{n})$, then $\phi\in
WPAA_{0}(\mathbb{T},\nu)$, where $\nu\in\mathcal{U}_{B}$ if and only
if, for every $\epsilon>0$,
\begin{equation*}
\lim\limits_{r\longrightarrow+\infty}\frac{1}{m(r,\nu,t_{0})}\nu_{\Delta}(M_{r,\epsilon,t_{0}}(\phi))=0,
\end{equation*}
where $r\in\Pi$ and
$M_{r,\epsilon,t_{0}}(\phi):=\left\{t\in[t_{0}-r,t_{0}+r]_{\mathbb{T}}:
\,\ \|\phi(t)\|\geq\epsilon\right\}$.
\end{lemma}
Proof. The demonstration is similar to the proof of Lemma 3.2 in
\cite{tomorphic2}.
\begin{lemma}
$WPAA_{0}(\mathbb{T},\nu)$ is a translation invariant set of
$S^{p}BC_{rd}(\mathbb{T},\mathbb{R}^{n})$ with respect to $\Pi$ if
$\nu\in\mathcal{U}_{B}$, i.e. for any $s\in\Pi$, one has
\begin{equation*}
\phi(t+s):=\theta_{s}\phi\in WPAA_{0}(\mathbb{T},\nu)\,\ \text{if}
\,\ \nu\in\mathcal{U}_{B}.
\end{equation*}
\end{lemma}
Proof. Similar to proof of Lemma 3.3 in \cite{tomorphic2}.
\begin{lemma}
Let $\phi\in S^{p}AA(\mathbb{T},\mathbb{R}^{n})$, then the range of
$\phi$, $\phi(\mathbb{T})$ is a relatively compact subset
$\mathbb{R}^{n}$.
\end{lemma}
Proof. Similar to proof of Lemma 3.3 in \cite{tomorphic2}.
\begin{lemma}\label{lem5}
If $f=g+\phi$ with $g\in S^{p}AA(\mathbb{T},\mathbb{R}^{n})$ and
$\phi\in WPAA_{0}(\mathbb{T},\nu)$, where $\nu\in\mathcal{U}_{B}$,
then $g(\mathbb{T})\subset \overline{f(\mathbb{T})}$.
\end{lemma}
Proof. The demonstration is similar to the proof of Lemma 3.5 in
\cite{tomorphic2}.
\begin{lemma}
The decomposition of a Stepanov-like weighted pseudo-almost
automorphic on time-space scales function according to
$S^{p}AA\oplus WPAA_{0}$ is unique for any $\nu\in\mathcal{U}_{B}$.
\end{lemma}
Proof. Assume that $f_{1}=g_{1}+\phi_{1}$ and
$f_{2}=g_{2}+\phi_{2}$. Then $(g_{1}-g_{2})+(\phi_{1}-\phi_{2})=0$.
Since $g_{1}-g_{2}\in S^{p}AA(\mathbb{T},\mathbb{R}^{n})$ , and
$\phi_{1}-\phi_{2}\in WPAA_{0}$ in view of Lemma \ref{lem5}, we
deduce that $g_{1}-g_{2}=0$. Consequently $\phi_{1}-\phi_{2}=0$,
i.e. $\phi_{1}=\phi_{2}$. This completes the proof.
\begin{lemma}
For $\nu\in\mathcal{U}_{B}$,
$(S^{p}WPAA(\mathbb{T},\nu),\|.\|_{S_{l}^{p}})$ is a Banach space.
\end{lemma}
Proof. Assume that $\{f_{n}\}_{n\in\mathbb{N}}$ is a Cauchy sequence
in $S^{p}WPAA(\mathbb{T},\nu)$. We can write uniquely
$f_{n}=g_{n}+\phi_{n}$. Using Lemma \ref{lem5}, we see that
$\|g_{p}-g_{q}\|\leq\|f_{p}-f_{q}\|_{\infty}$, from which we deduce
that $\{g_{n}\}_{n\in\mathbb{N}}$ is a Cauchy sequence in
$AA(\mathbb{T},\mathbb{R}^{n})$. Hence, $\phi_{n}=f_{n}-g_{n}$ is a
Cauchy sequence in $WPAA_{0}(\mathbb{T},\nu)$. We deduce that
$g_{n}\longrightarrow g\in AA(\mathbb{T},\mathbb{R}^{n})$,
$\phi_{n}\longrightarrow \phi\in WPAA(\mathbb{T},\nu)$ and finally
$f_{n}\longrightarrow g+\phi\in S^{p}WPAA(\mathbb{T},\nu)$. This
complete the proof.
\begin{definition}(\cite{tomorphic2})
Let $\nu_{1},\nu_{2}\in\mathcal{U}_{\infty}$. One says that
$\nu_{1}$ is equivalent to $\nu_{2}$, written $\nu_{1}\sim \nu_{2}$
if $\nu_{1}/\nu_{2}\in \mathcal{U}_{B}$.
\end{definition}
\begin{lemma}
Let $\nu_{1},\nu_{2}\in\mathcal{U}_{\infty}$. If $\nu_{1}\sim
\nu_{2}$, then
$S^{p}WPAA(\mathbb{T},\nu_{1})=S^{p}WPAA(\mathbb{T},\nu_{2})$.
\end{lemma}
Proof. The demonstration is similar to the proof of Theorem 3.8 in
\cite{tomorphic2}.
\begin{lemma}
Let $f=g+\phi\in S^{p}WPAA(\mathbb{T},\nu)$, where
$\nu\in\mathcal{U}_{B}$. Assume that $f$ and $g$ are Lipshitzian in
$x\in\mathbb{R}^{n}$ uniformly in $t\in\mathbb{T}$, then
$f(.,h(.))\in S^{p}WPAA(\mathbb{T},\nu)$ if $h\in
S^{p}WPAA(\mathbb{T},\nu)$.
\end{lemma}
Proof. Similar to proof of Theorem 3.10 in \cite{tomorphic2}.
\begin{lemma}\label{uniform} (\cite{diaganna})
If $f(t)$ is almost automorphic, $F(.)$ is uniformly continuous on
the value field of $f(t)$, then $F\circ f$ is almost automorphic.
\end{lemma}
\begin{lemma}
If $f\in S^{p}C(\mathbb{R},\mathbb{R})$ satisfies the Lipschitz
condition (with $L$ is a constant of Lipschitz), $\varphi\in
S^{p}WPAA(\mathbb{T},\nu)$, $\theta\in
S^{p}C_{rd}^{1}(\mathbb{T},\Pi)$ and
$\eta:=\inf\limits_{t\in\mathbb{T}}\left(1-\theta^{\Delta}(t)\right)>0$,
then $f(\varphi(t-\theta(t)))\in S^{p}WPAA(\mathbb{T},\nu)$.
\end{lemma}
Proof.\\
From Definition \ref{defwpap}, we have
$\varphi=\varphi_{1}+\varphi_{2}$, where $\varphi_{1}\in
AP(\mathbb{T})$ and $\varphi_{2}\in WPAP_{0}(\mathbb{T},\nu,t_{0})$.
Set
\begin{eqnarray*}
E(t)&=& f(\varphi(t-\theta(t)))=f(\varphi_{1}(t-\theta(t)))+\left[f(\varphi_{1}(t-\theta(t))-\varphi_{2}(t-\theta(t)))\right.\\
&-&\left.f(t-\varphi_{1}(t-\theta(t)))\right]=E_{1}(t)+E_{2}(t).
\end{eqnarray*}
Firstly, it follows from Lemma \ref{uniform} that $E_{1}\in
S^{p}AP(\mathbb{T})$. Next, we show that $E_{2}\in
WPAA_{0}(\mathbb{T},\nu,t_{0})$. Since
\begin{eqnarray*}
&&\lim\limits_{r\longrightarrow+\infty}\frac{1}{m(r,\nu,t_{0})}\int_{t_{0}-r}^{t_{0}+r}\left\vert E_{2}(s)\right\vert \nu(s)\Delta s\\
&=& \lim\limits_{r\longrightarrow+\infty}\frac{1}{m(r,\nu,t_{0})}\int_{t_{0}-r}^{t_{0}+r}\left\vert f(\varphi_{1}(t-\theta(t))-\varphi_{2}(t-\theta(t)))\right.\\
&-&\left.f(t-\varphi_{1}(t-\theta(t)))\right\vert\nu(s)\Delta s\\
&\leq&
\lim\limits_{r\longrightarrow+\infty}\frac{L}{m(r,\nu,t_{0})}\int_{t_{0}-r}^{t_{0}+r}\left\vert\varphi_{2}(t-\theta(t)))\right\vert
\nu(s)\Delta s
\end{eqnarray*}
and
\begin{eqnarray*}
0&\leq&\frac{L}{m(r,\nu,t_{0})}\int_{t_{0}-r}^{t_{0}+r}\left\vert\varphi_{2}(t-\theta(t)))\right\vert \nu(s)\Delta s\\
&=& \frac{L}{m(r,\nu,t_{0})}\int_{t_{0}-r-\theta(t_{0}-r)}^{t_{0}+r-\theta(t_{0}+r)}\frac{1}{1-\theta^{\Delta}(s)}\left\vert\varphi_{2}(u))\right\vert \nu(u)\Delta u\\
&\leq&
\frac{1}{\eta}\frac{r+\theta^{+}}{r}\frac{L}{m(r,\nu,t_{0})}\int_{t_{0}-(r+\theta^{+})}^{t_{0}+r+\theta^{+}}\left\vert\varphi_{2}(u))\right\vert
\nu(u)\Delta u=0,
\end{eqnarray*}
$E_{2}\in WPAA_{0}(\mathbb{T},\nu,t_{0})$. Thus $E\in
S^{p}WPAA(\mathbb{T},\nu)$. The proof is achieved.
\section{Model description and hypotheses}
In this paper, we consider a class of $n$-neuron high-order BAM
neural networks (HOBAMs) with mixed time-varying delays and leakage
delays on time-space scales which are defined in the following
lines:
\begin{equation}\label{eq1}
\left\{
\begin{array}{ccc}
x^{\Delta}_{i}\left( t\right) =
-\alpha_{i}(t)x_{i}\left(t-\eta_{i}(t)\right)
+\sum\limits_{j=1}^{m}D_{ij}\left( t\right) f_{j}\left(y_{j}\left(
t\right) \right) +\sum\limits_{j=1}^{m}D_{ij}^{\tau}\left( t\right)
f_{j}\left(x_{j}\left( t -\tau_{ij}(t)\right) \right)\\
+\sum\limits_{j=1}^{m}\overline{D}_{ij}\left(
t\right)\int_{t-\sigma_{ij}(t)} ^{t}f_{j}\left(y_{j}\left( s\right)
\right)\Delta s +\sum\limits_{j=1}^{m}\widetilde{D}_{ij}\left(
t\right)\int_{t-\xi_{ij}(t)} ^{t}f_{j}\left(y_{j}^{\Delta}\left(
s\right) \right)\Delta s\\
+\sum\limits_{j=1}^{m}\sum\limits_{k=1}^{m}T_{ijk}(t)f_{k}(y_{k}(t-\chi_{k}(t)))f_{j}(y_{j}(t-\chi_{j}(t)))
+I_{i}\left( t\right),\\
y^{\Delta}_{j}\left( t\right) =
-c_{j}(t)y_{j}\left(t-\eta_{j}(t)\right)
+\sum\limits_{i=1}^{n}E_{ij}\left( t\right) f_{j}\left(x_{j}\left(
t\right) \right) +\sum\limits_{i=1}^{n}E_{ij}^{\tau}\left( t\right)
f_{j}\left(x_{j}\left( t -\tau_{ij}(t)\right) \right)\\
+\sum\limits_{i=1}^{n}\overline{E}_{ij}\left(
t\right)\int_{t-\sigma_{ij}(t)} ^{t}f_{j}\left(x_{j}\left( s\right)
\right)\Delta s +\sum\limits_{i=1}^{n}\widetilde{E}_{ij}\left(
t\right)\int_{t-\xi_{ij}(t)}^{t}f_{j}\left(x_{j}^{\Delta}\left(
s\right) \right)\Delta s\\
+\sum\limits_{i=1}^{n}\sum\limits_{k=1}^{n}\overline{T}_{ijk}(t)f_{k}(x_{k}(t-\chi_{k}(t)))f_{j}(x_{j}(t-\chi_{j}(t)))
+J_{j}\left( t\right),\,\ t\in\mathbb{T},
\end{array}
\right.
\end{equation}
where $i=1,\cdots ,n$ and $j=1,\cdots ,m$; $\mathbb{T}$ is an almost
periodic time scale; $x_{i}(t)$ and $y_{i}(t)$ are the neuron
current activity level of ith neuron in the first layer and the jth
neuron in the second layer respectively at time $t$ ($i=1,\cdots ,n,
\,\ j=1,\cdots ,m$); $\alpha_{i}(t), c_{j}(t)$ are the time variable
of the neuron $i$ in the first layer and the neuron $j$ in the
second neuron respectively; $f_{i}(x_{i}(t))$ and $f_{j}(x_{j}(t))$
are the output of neurons; $I_{i}(t)$ and $J_{j}(t)$ denote the
external inputs on the ith neuron at time $t$ for the first layer
and the jth neuron at the second layer at time $t$;
\begin{multline*} \left.
\begin{array}{ll}
t\longmapsto D_{ij}(t)\\
t\longmapsto D_{ij}^{\tau}(t)\\
t\longmapsto \overline{D}_{ij}(t)\\
t\longmapsto \widetilde{D}_{ij}(t)\\
t\longmapsto E_{ij}(t)\\
t\longmapsto E_{ij}^{\tau}(t)\\
t\longmapsto \overline{E}_{ij}(t)\\
t\longmapsto \widetilde{E}_{ij}(t)
\end{array}\right\} \text{represent the connection weights and the synaptic weights }\\
\text{of delayed feedback between the $i$th neuron and the $j$th
neuron respectively};
\end{multline*}
for all $i,k=1,...n, \,\ j=1,...,m \,\ t\longmapsto T_{ijk}(t)$ and $t\longmapsto\overline{T}_{ijk}(t)$ are the second-order connection weights of delayed feedback;\\
$t\longmapsto I_{i}(t)$, $t\longmapsto J_{i}(t)$ denote the external
inputs on the $i$th neuron at time $t$; $t\longmapsto \eta_{i}(t)$
and $t\longmapsto \varsigma_{i}(t)$ are leakage delays and satisfy
$t-\eta_{i}(t)\in\mathbb{T}$, $t-\varsigma_{i}(t)\in\mathbb{T}$ for
$t\in\mathbb{T}$; $\tau_{ij}(t)$, $\sigma_{ij}(t)$, $\chi_{k}(t)$
and $\xi_{ij}(t)$ are transmission delays and satisfy
$t-\tau_{ij}(t)\in\mathbb{T}$, $t-\sigma_{ij}(t)\in\mathbb{T}$,
$t-\chi_{k}(t)\in\mathbb{T}$ and $t-\xi_{ij}(t)\in\mathbb{T}$ for
$t\in\mathbb{T}$.

For convenience, we introduce the following notations:
\begin{multline*}
\alpha_{i}^{+}=\sup\limits_{t\in\mathbb{T}}|\alpha_{i}(t)|, \,\ \alpha_{i}^{-}=\inf\limits_{t\in\mathbb{T}}|\alpha_{i}(t)|>0,\,\ c_{i}^{+}=\sup\limits_{t\in\mathbb{T}}|c_{i}(t)|,\,\ c_{i}^{-}=\inf\limits_{t\in\mathbb{T}}|c_{i}(t)|>0,\\
\eta_{i}^{+}=\sup\limits_{t\in\mathbb{T}}|\eta_{i}(t)|, \,\
\varsigma_{i}^{+}=\sup\limits_{t\in\mathbb{T}}|\varsigma_{i}(t)|,\,\
D_{i}^{+}=\sup\limits_{t\in\mathbb{T}}|D_{i}(t)|,\,\
(D_{i}^{\tau})^{+}=\sup\limits_{t\in\mathbb{T}}|D_{i}^{\tau}(t)|,\\
D_{ij}^{+}=\sup\limits_{t\in\mathbb{T}}|D_{ij}(t)|,\,\
(D_{ij}^{\tau})^{+}=\sup\limits_{t\in\mathbb{T}}|D_{ij}^{\tau}(t)|,\,\
\overline{D}_{ij}^{+}=\sup\limits_{t\in\mathbb{T}}|\overline{D}_{ij}(t)|,\\
(\widetilde{D}_{ij})^{+}=\sup\limits_{t\in\mathbb{T}}|\widetilde{D}_{ij}(t)|,\,\
E_{i}^{+}=\sup\limits_{t\in\mathbb{T}}|E_{i}(t)|,\,\
(E_{i}^{\tau})^{+}=\sup\limits_{t\in\mathbb{T}}|E_{i}^{\tau}(t)|,\,\
E_{ij}^{+}=\sup\limits_{t\in\mathbb{T}}|E_{ij}(t)|,\\
(E_{ij}^{\tau})^{+}=\sup\limits_{t\in\mathbb{T}}|E_{ij}^{\tau}(t)|,\,\
\overline{E}_{ij}^{+}=\sup\limits_{t\in\mathbb{T}}|\overline{E}_{ij}(t)|,\,\
(\widetilde{E}_{ij})^{+}=\sup\limits_{t\in\mathbb{T}}|\widetilde{E}_{ij}(t)|,\\
T_{ijk}^{+}=\sup\limits_{t\in\mathbb{T}}|T_{ijk}(t)|,\,\
\overline{T}_{ijk}^{+}=\sup\limits_{t\in\mathbb{T}}|\overline{T}_{ijk}(t)|,\,\
\tau_{ij}^{+}=\sup\limits_{t\in\mathbb{T}}|\tau_{ij}(t)|,\\
\sigma_{ij}^{+}=\sup\limits_{t\in\mathbb{T}}|\sigma_{ij}(t)|,\,\
\xi_{ij}^{+}=\sup\limits_{t\in\mathbb{T}}|\xi_{ij}(t)|,\,\
\chi_{j}^{+}=\sup\limits_{t\in\mathbb{T}}|\chi_{j}(t)|,\,\
i,=1,...,n, \,\ j=1,...,m.
\end{multline*}

We denote that $[a,b]_{\mathbb{T}}=\{t, t\in[a,b]\cap \mathbb{T}\}$.
The initial conditions associated with system (\ref{eq1}), are of
the form:
\begin{equation*}
x_{i}\left( s\right) = \varphi_{i}\left(s\right),\,\
y_{j}\left(s\right)= \phi_{j}\left( s\right),\,\ s\in [-\theta
,0]_{\mathbb{T}} ,1\leq i\leq n,\,\ 1\leq j\leq m,
\end{equation*}
where $\varphi_{i}(.)$ and $\phi_{i}(.)$ are the real-valued bounded
$\Delta$-differentiable functions defined on
$[-\theta,0]_{\mathbb{T}}$,
\begin{multline*}
\theta=\max\{\eta,\tau,\chi,\sigma,\xi,\varsigma\}, \,\
\eta=\max\limits_{1\leq i\leq n}\eta_{i}^{+}, \,\
\tau=\max\limits_{1\leq i\leq n, 1\leq j\leq m}\tau_{ij}^{+}, \,\
\chi=\max\limits_{1\leq j\leq m}\chi_{j}^{+},\\
\sigma=\max\limits_{1\leq i\leq n, 1\leq j\leq m}\sigma_{ij}^{+},\,\
\xi=\max\limits_{1\leq i\leq n, 1\leq j\leq m}\xi_{ij}^{+}\,\
\text{and} \,\ \varsigma=\max\limits_{1\leq i\leq n, 1\leq j\leq
m}\varsigma_{ij}^{+}.
\end{multline*}
\begin{remark}
This is the first time to study the Stepanov-like weighted
pseudo-almost automorphic solutions of system (\ref{eq1}) for the
both cases: continuous and discrete. Furthermore, there is no result
about automorphic, Stepanov almost automorphic and Stepanov-like
weighted pseudo-almost automorphic solutions of networks
(\ref{eq1}).
\end{remark}
Let us list some assumptions that will be used throughout the rest
of this paper.
\begin{description}
\item[($H_{1}$)]
For all $1\leq i,j\leq n,$, the functions $\alpha_{i}\left(\cdot
\right), c_{j}\left(\cdot \right)\in\mathcal{R}_{\nu}^{+}$ and
$D_{ij}(.)$, $D^{\tau}_{ij}(.)$ $\overline{D}_{ij}(.)$,
$\widetilde{D}_{ij}(.)$, $T_{ijk}(.)$, $E_{ij}(.)$,
$E^{\tau}_{ij}(.)$ $\overline{E}_{ij}(.)$, $\widetilde{E}_{ij}(.)$,
$\overline{T}_{ijk}(.)$, $\eta_{i}(.)$, $\varsigma_{j}(.)$,
$\tau_{ij}(.)$, $\chi_{j}(.)$, $\sigma_{ij}(.)$, $\xi_{ij}(.)$,
$I_{i}(.)$, $J_{j}(.)$ are $ld$-continuous Stepanov-like weighted
pseudo-almost automorphic
functions for $i=1,...,n$, $j=1,...,m$.\\
\item[($H_{2}$)]
The functions $f_{j}\left(\cdot\right)$ are $\Delta$-differential
and satisfy the Lipschitz condition, i.e., there are constants
$L_{j}>0$ such that for all $x,y\in\mathbb{R}$, and for all $1\leq
j\leq \max\{n,m\},$ one has $ \left\vert f_{j}\left( x\right)
-f_{j}\left( y\right) \right\vert \leq L_{j}\left\vert
x-y\right\vert. $
\item[($H_{3}$)]
\begin{multline*}
\phantom{++++++}\max\limits_{1\leq i\leq
n}\left\{\frac{M_{i}}{\alpha_{i}^{-}},\left(1+\frac{\alpha_{i}^{+}}{\alpha_{i}^{-}}\right)M_{i},\frac{N_{i}}{c_{i}^{-}},
\left(1+\frac{c_{i}^{+}}{c_{i}^{-}}\right)N_{i}\right\}\leq r \\
\text{and} \,\ \max\limits_{1\leq i\leq
n}\left\{\frac{\overline{M}_{i}}{\alpha_{i}^{-}},\left(1+\frac{\alpha_{i}^{+}}{\alpha_{i}^{-}}\right)\overline{M}_{i},
\frac{\overline{N}_{i}}{c_{i}^{-}},\left(1+\frac{c_{i}^{+}}{c_{i}^{-}}\right)\overline{N}_{i}\right\}\leq
1,\phantom{++}
\end{multline*}
where $r$ is a constant, for $i=1,...,n$ and $j=1,...,m$,
\begin{eqnarray*}
M_{i}&=&\alpha_{i}^{+}\eta_{i}^{+}r+\sum\limits_{j=1}^{m}\left(D_{ij}^{+}+(D_{ij}^{\tau})^{+}+\overline{D}_{ij}^{+}\sigma_{ij}^{+}
+\widetilde{D}_{ij}^{+}\xi_{ij}^{+}\right)(L_{j}r+|f_{j}(0)|)\\
&+&\sum\limits_{j=1}^{m}\sum\limits_{k=1}^{m}T_{ijk}^{+}\left(L_{k}r+\left\vert f_{k}(0)\right\vert\right)\left(L_{j}r+\left\vert f_{j}(0)\right\vert\right)+I_{i}^{+},\\
\overline{M}_{i}&=&\alpha_{i}^{+}\eta_{i}^{+}+\sum\limits_{j=1}^{m}\left(D_{ij}^{+}+(D_{ij}^{\tau})^{+}
+\overline{D}_{ij}^{+}\sigma_{ij}^{+}+\widetilde{D}_{ij}^{+}\xi_{ij}^{+}\right)L_{j}\\
&+&\sum\limits_{j=1}^{m}\sum\limits_{k=1}^{m}(T_{ijk}^{+}+T_{ikj}^{+})(L_{k}r+|f_{k}(0)|),\\
N_{j}&=&c_{j}^{+}\varsigma_{j}^{+}r+\sum\limits_{i=1}^{n}\left(E_{ij}^{+}+(E_{ij}^{\tau})^{+}+\overline{E}_{ij}^{+}\sigma_{ij}^{+}
+\widetilde{E}_{ij}^{+}\xi_{ij}^{+}\right)(L_{i}r+|f_{i}(0)|)\\
&+&\sum\limits_{i=1}^{n}\sum\limits_{k=1}^{n}\overline{T}_{ijk}^{+}\left(L_{k}r+\left\vert f_{k}(0)\right\vert\right)\left(L_{i}r+\left\vert f_{i}(0)\right\vert\right)+J_{j}^{+},\\
\overline{N}_{j}&=&
c_{j}^{+}\varsigma_{j}^{+}+\sum\limits_{i=1}^{n}\left(E_{ij}^{+}+(E_{ij}^{\tau})^{+}
+\overline{E}_{ij}^{+}\sigma_{ij}^{+}+\widetilde{E}_{ij}^{+}\xi_{ij}^{+}\right)L_{i}\\
&+&\sum\limits_{i=1}^{n}\sum\limits_{k=1}^{n}(\overline{T}_{ijk}^{+}+\overline{T}_{ikj}^{+})(L_{k}r+|f_{k}(0)|).
\end{eqnarray*}
\item[($H_{4}$)]
$\inf\limits_{t\in\mathbb{T}}\left(1-\sigma_{ij}^{\nabla}(t)\right)>0$,
$\inf\limits_{t\in\mathbb{T}}\left(1-\xi_{ij}^{\nabla}(t)\right)>0$,
and for all $s\in\Pi$,
\begin{equation*}
\limsup\limits_{|t|\longrightarrow+\infty}\frac{\nu(t+s)}{\nu(t)}<\infty.
\end{equation*}
\end{description}
\begin{remark}
The bidirectional associative memory (BAM) neural networks with
mixed time-varying delays and leakage time-varying delays on
time-space scales is investigated in \cite{maas}. Some sufficient
conditions are given for the existence, convergence and the global
exponential stability of the weighted pseudo almost-periodic
solution. However, Theorem 4.1, Theorem 5.1 and Theorem 6.1 proposed
in \cite{maas} are not applicable for the HOBAMs with mixed
time-varying delays in the leakage terms.
\end{remark}
\section{The existence of Stepanov-like weighted pseudo-almost automorphic on time-space scales solutions}
In this section, based on Banach's fixed point theorem and the
theory of calculus on time-space scales, we will present a new
condition for the existence and uniqueness of weighted pseudo-almost
automorphic on time-space scales solutions of (\ref{eq1}).
Additionally, we will show a result about the delta derivative of
the only Stepanov-like weighted pseudo-almost automorphic on
time-space scales solution of system (\ref{eq1}).

Let
\begin{equation*}
\mathbb{B}=\{(\varphi_{1},\varphi_{2},...,\varphi_{n},\phi_{1},\phi_{2},...,\phi_{m})^{T}:\varphi_{i},\phi_{j}\in
C^{1}(\mathbb{T},\mathbb{R}),\,\ i=1,...,n, \,\ j=1,...,m\}.
\end{equation*}
with the norm $
\|\psi\|_{\mathbb{B}}=\sup\limits_{t\in\mathbb{T}}\max\limits_{i=1,...,n\,\
j=1,...,m}\{|\varphi_{i}(t)|,|\phi_{j}(t)|,|\varphi^{\Delta}_{i}(t)|,|\phi^{\Delta}_{j}(t)|\},
$
then $(\mathbb{B},\|\psi\|_{\mathbb{B}})$ is a Banach space.\\

For every
$\psi=(\varphi_{1},...\varphi_{n},\phi_{1},...,\phi_{m})\in\mathbb{B}$,
we consider the following system
\begin{equation}
x_{i}^{\Delta}(t)=-\alpha_{i}(t)x_{i}(t)+F_{i}(t,\varphi_{i}),\,\
y_{j}^{\Delta}(t)=-c_{j}(t)y_{j}(t)+G_{j}(t,\phi_{i}), \,\
t\in\mathbb{T},
\end{equation}
where, for $i=1,...,n$ and $j=1,...,m$
\begin{eqnarray*}
F_{i}(t,\varphi_{i}(t))&=&\alpha_{i}(t)\int_{t-\eta_{i}(t)}^{t}\varphi_{i}^{\Delta}(s)\Delta s +\sum\limits_{j=1}^{m}D_{ij}\left( t\right) f_{j}\left(\phi_{j}\left( t\right) \right)\\
&+&\sum\limits_{j=1}^{m}D_{ij}^{\tau}\left(t\right)
f_{j}\left(\phi_{j}\left( t -\tau_{ij}(t)\right) \right)
+\sum\limits_{j=1}^{m}\overline{D}_{ij}\left(t\right)\int_{t-\sigma_{ij}(t)}^{t}f_{j}\left(\phi_{j}\left( s\right) \right)\Delta s\\
&+&\sum\limits_{j=1}^{m}\widetilde{D}_{ij}\left( t\right)\int_{t-\xi_{ij}(t)} ^{t}f_{j}\left(\phi_{j}^{\Delta}\left( s\right) \right)\Delta s\\
&+&\sum\limits_{j=1}^{m}\sum\limits_{k=1}^{n}T_{ijk}(t)f_{k}(\phi_{k}(t-\chi_{k}(t)))f_{j}(\phi_{j}(t-\chi_{j}(t)))+I_{i}(t),\\
G_{j}(t,\phi_{j}(t))&=&c_{j}(t)\int_{t-\varsigma_{j}(t)}^{t}\phi_{j}^{\Delta}(s)\Delta
s+\sum\limits_{i=1}^{n}E_{ij}\left( t\right) f_{i}\left(\varphi_{i}\left( t\right) \right)\\
&+&\sum\limits_{i=1}^{n}E_{ij}^{\tau}\left(t\right)
f_{i}\left(\varphi_{i}\left( t -\tau_{ij}(t)\right) \right)
+\sum\limits_{i=1}^{n}\overline{E}_{ij}\left(t\right)\int_{t-\sigma_{ij}(t)}^{t}f_{i}\left(\varphi_{i}\left( s\right) \right)\Delta s\\
&+&\sum\limits_{i=1}^{n}\widetilde{E}_{ij}\left( t\right)\int_{t-\xi_{ij}(t)} ^{t}f_{i}\left(\varphi_{i}^{\Delta}\left( s\right) \right)\Delta s\\
&+&\sum\limits_{i=1}^{n}\sum\limits_{k=1}^{n}\overline{T}_{ijk}(t)f_{k}(\varphi_{k}(t-\chi_{k}(t)))f_{i}(\varphi_{i}(t-\chi_{i}(t)))+J_{j}(t).
\end{eqnarray*}
Let
$y_{\psi}(t)=\left(x_{\varphi_{1}}(t),...,x_{\varphi_{n}}(t),y_{\phi_{1}}(t),...,y_{\phi_{m}}(t)\right)^{T}$,
where:
\begin{eqnarray*}
x_{\varphi_{i}}(t)&=&\int_{-\infty}^{t}\hat{e}_{-\alpha_{i}}(t,\sigma(s))F_{i}(t,\varphi_{i}(s))\Delta s,\\
y_{\phi_{j}}(t)&=&\int_{-\infty}^{t}\hat{e}_{-c_{j}}(t,\sigma(s))G_{j}(t,\phi_{j}(s))\Delta
s,
\end{eqnarray*}
\begin{lemma}
Suppose that assumptions $(H_{1})-(H_{4})$ hold. Define the
nonlinear operator $\Gamma:\mathbb{F}\longrightarrow \mathbb{F}$ by
for each $\psi\in WPAA(\mathbb{T},\nu)$
\begin{equation*}
(\Gamma\psi)(t)=y_{\psi}(t), \,\ \psi\in\mathbb{F}.
\end{equation*}
Then $\Gamma$ maps $WPAA(\mathbb{T},\nu)$ into itself.
\end{lemma}
Proof.\\
We show that for any $\psi\in\mathbb{F}$, $\Gamma\psi\in\mathbb{F}$.
\begin{eqnarray*}
\left\vert F_{i}(s,\varphi_{i}(s))\right\vert
&\leq&\alpha_{i}^{+}\eta_{i}^{+}r
+\sum\limits_{j=1}^{n}D_{ij}^{+}\left(L_{j}\left\vert\varphi_{j}\left(s\right)\right\vert+\left\vert f_{j}(0)\right\vert\right)\\
&+&\sum\limits_{j=1}^{n}(D_{ij}^{\tau})^{+}\left(L_{j}\left\vert\varphi_{j}\left(s-\tau_{ij}(s)\right)\right\vert+\left\vert f_{j}(0)\right\vert\right)\\
&+&\sum\limits_{j=1}^{n}\overline{D}_{ij}^{+}\sigma_{ij}^{+}\left(L_{j}r+\left\vert
f_{j}(0)\right\vert\right)+
\sum\limits_{j=1}^{n}\widetilde{D}_{ij}^{+}\xi_{ij}^{+}\left(L_{j}r+\left\vert f_{j}(0)\right\vert\right)\\
&+&
\sum\limits_{j=1}^{n}\sum\limits_{k=1}^{n}T_{ijk}^{+}\left(L_{k}r+\left\vert
f_{k}(0)\right\vert\right)\left(L_{j}r+\left\vert f_{j}(0)\right\vert\right)+I_{i}^{+}\\
&\leq& M_{i}.
\end{eqnarray*}
In a similar way, we have
\begin{equation*}
\left\vert G_{j}(s,\phi_{j}(s))\right\vert\leq N_{j}.
\end{equation*}
Which leads to, for $i=1,...,n$, $j=1,...,m$.
\begin{eqnarray*}
\sup\limits_{t\in\mathbb{T}}\left\vert x_{\varphi_{i}}(t)\right\vert&=&\sup\limits_{t\in\mathbb{T}}\left\vert\int_{-\infty}^{t}\hat{e}_{-\alpha_{i}^{-}}(t,\sigma(s))F_{i}(t,\varphi_{i}(s))\Delta s\right\vert\\
&\leq& \sup\limits_{t\in\mathbb{T}}\int_{-\infty}^{t}\hat{e}_{-\alpha_{i}^{-}}(t,\sigma(s))\left\vert F_{i}(s,\varphi_{i}(s))\right\vert \Delta s\\
&\leq& \frac{M_{i}}{\alpha_{i}^{-}},
\end{eqnarray*}
and
\begin{equation*}
\sup\limits_{t\in\mathbb{T}}\left\vert
y_{\phi_{j}}(t)\right\vert=\sup\limits_{t\in\mathbb{T}}\left\vert\int_{-\infty}^{t}\hat{e}_{-c_{i}}(t,\sigma(s))G_{i}(t,\phi_{i}(s))\Delta
s\right\vert\leq \frac{N_{j}}{c_{j}^{-}}.
\end{equation*}
Otherwise, for $i=1,...,n$, we have
\begin{eqnarray*}
\sup\limits_{t\in\mathbb{T}}\left\vert
x_{\varphi_{i}}^{\Delta}(t)\right\vert
&=&\sup\limits_{t\in\mathbb{T}}\left\vert\left(\int_{-\infty}^{t}\hat{e}_{-\alpha_{i}}(t,\sigma(s))F_{i}(t,\varphi_{i}(s))\Delta s\right)^{\Delta} \right\vert\\
&=&\sup\limits_{t\in\mathbb{T}}\left\vert F_{i}(t,\varphi_{i}(t))-\alpha_{i}(t)\int_{-\infty}^{t}\hat{e}_{-\alpha_{i}}(t,\sigma(s))F_{i}(t,\varphi_{i}(s))\Delta s\right\vert\\
&\leq&\alpha_{i}^{+}\eta_{i}^{+}r+\sum\limits_{j=1}^{m}\left(D_{ij}^{+}+(D_{ij}^{\tau})^{+}+(\overline{D}_{ij})^{+}\sigma_{ij}^{+}
+(\widetilde{D}_{ij})^{+}\xi_{ij}^{+}\right)(L_{j}r+\left\vert f_{j}(0)\right\vert)\\
&+&\sum\limits_{j=1}^{m}\sum\limits_{k=1}^{n}T_{ijk}^{+}\left(L_{k}r+\left\vert f_{k}(0)\right\vert\right)\left(L_{j}r+\left\vert f_{j}(0)\right\vert\right)\\
&+&I_{i}^{+}+\frac{\alpha_{i}^{+}}{\alpha_{i}^{-}}\left(\alpha_{i}^{+}\eta_{i}^{+}r
+\sum\limits_{j=1}^{m}\left(D_{ij}^{+}+(D_{ij}^{\tau})^{+}+(\overline{D}_{ij})^{+}\sigma_{ij}^{+}\right.\right.\\
&+&\left.(\widetilde{D}_{ij})^{+}\xi_{ij}^{+}\right)(L_{j}r+\left\vert
f_{j}(0)\right\vert)+\sum\limits_{j=1}^{m}\sum\limits_{k=1}^{n}T_{ijk}^{+}\left(L_{k}r+\left\vert
f_{k}(0)\right\vert\right)\\
&\times&\left(L_{j}r+\left\vert f_{j}(0)\right\vert\right)
+\left. I_{i}^{+}\right)\\
&=&\left(1+\frac{\alpha_{i}^{+}}{\alpha_{i}^{-}}\right)M_{i}.
\end{eqnarray*}
Similarly,
\begin{equation*}
\sup\limits_{t\in\mathbb{T}}\left\vert
y_{\phi_{j}}^{\Delta}(t)\right\vert\leq
\left(1+\frac{c_{j}^{+}}{c_{j}^{-}}\right)N_{j}.
\end{equation*}
From hypothesis $(H_{3})$, we can obtain
\begin{equation*}
\|\Gamma \psi\|_{\mathbb{B}}\leq r,
\end{equation*}
which implies that operator $\Gamma$ is a self-mapping from $\mathbb{F}$ to $\mathbb{F}$.\\
\begin{theorem}\label{th0}
Let $(H_{1})-(H_{4})$ hold. The system (\ref{eq1}) has a unique
Stepanov-like weighted pseudo-almost automorphic solution in the
region $\mathbb{F}=\{\psi\in\mathbb{B}: \|\psi\|_{\mathbb{B}}\leq
r\}.$
\end{theorem}
Proof.\\
First, for
$\psi=\left(\varphi_{1},...,\varphi_{n},\phi_{1},...,\phi_{m}\right)^{T}$,
$\Omega=\left(u_{1},...,u_{n},
v_{1},...,v_{m}\right)^{T}\in\mathbb{F}$, we have
\begin{eqnarray*}
\sup\limits_{s\in\mathbb{T}}\left\Vert
x_{\varphi_{i}}(s)-x_{u_{i}}(s) \right\Vert
&\leq&\frac{1}{\alpha_{i}^{-}}\left(\alpha_{i}^{+}\eta_{i}^{+}+\sum\limits_{j=1}^{m}\left(D_{ij}^{+}+(D_{ij}^{\tau})^{+}
+(\overline{D}_{ij})^{+}\sigma_{ij}^{+}+(\widetilde{D}_{ij})^{+}\xi_{ij}^{+}\right)L_{j}\right.\\
&+&\left.\sum\limits_{j=1}^{m}\sum\limits_{k=1}^{m}(T_{ijk}^{+}+T_{ikj}^{+})(L_{k}r+|f_{k}(0)|)\right)\|\psi-\Omega\|_{\mathbb{B}}\\
&=&
\frac{\overline{M}_{i}}{\alpha_{i}^{-}}\|\psi-\Omega\|_{\mathbb{B}}.
\end{eqnarray*}
Besides,
\begin{eqnarray*}
\sup\limits_{s\in\mathbb{T}}\left\Vert
x^{\Delta}_{\varphi_{i}}(s)-x^{\Delta}_{u_{i}}(s) \right\Vert
&\leq&\left(\alpha_{i}^{+}\eta_{i}^{+}+\sum\limits_{j=1}^{m}\left(D_{ij}^{+}+(D_{ij}^{\tau})^{+}
+(\overline{D}_{ij})^{+}\sigma_{ij}^{+}+(\widetilde{D}_{ij})^{+}\xi_{ij}^{+}\right)L_{j}\right.\\
&+&\left.\sum\limits_{j=1}^{m}\sum\limits_{k=1}^{m}(T_{ijk}^{+}+T_{ikj}^{+})(L_{k}r+|f_{k}(0)|)\right)\|\psi-\Omega\|_{\mathbb{B}}\\
&+&\frac{\alpha_{i}^{+}}{\alpha_{i}^{-}}\left(\alpha_{i}^{+}\eta_{i}^{+}\left\Vert
\phi_{i}^{\Delta}(s)-u_{i}^{\Delta}(s)\right\Vert+\sum\limits_{j=1}^{m}
D_{ij}^{+}\left\Vert \phi_{i}(t)-u_{i}(t)\right\Vert\right.\\
&+&\sum\limits_{j=1}^{m}(D_{ij}^{\tau})^{+}L_{j}\left\Vert \phi_{j}(s-\tau_{ij}(s))-u_{j}(s-\tau_{ij}(s))\right\Vert\\
&+&\sum\limits_{j=1}^{m}(\overline{D}_{ij})^{+}L_{j}\sigma_{ij}^{+}\left\Vert \phi_{j}(s)-u_{j}(s)\right\Vert\\
&+&\sum\limits_{j=1}^{m}\sum\limits_{k=1}^{m}(T_{ijk}^{+}+T_{ikj}^{+})(L_{k}r+|f_{k}(0)|)\left\Vert\phi_{j}(s)-u_{j}(s)\right\Vert\\
&+&\left.\sum\limits_{j=1}^{m}(\widetilde{D}_{ij})^{+}L_{j}\xi_{ij}^{+}\left\Vert
\phi_{j}^{\Delta}(s)-u^{\Delta}_{j}(s)\right\Vert\right)\\
&\leq&
\left(1+\frac{\alpha_{i}^{+}}{\alpha_{i}^{-}}\right)\overline{M}_{i}\|\psi-\Omega\|_{\mathbb{B}}.
\end{eqnarray*}
Similarly,
\begin{equation*}
\sup\limits_{s\in\mathbb{T}}\left\Vert y_{\phi_{j}}(s)-y_{v_{j}}(s)
\right\Vert \leq
\left(\frac{\overline{N}_{j}}{c_{j}^{-}}\right)\|\psi-\Omega\|_{\mathbb{B}},
\end{equation*}
and
\begin{equation*}
\sup\limits_{s\in\mathbb{T}}\left\Vert
y^{\Delta}_{\phi_{j}}(s)-y^{\Delta}_{v_{j}}(s) \right\Vert \leq
\left(1+\frac{c_{j}^{+}}{c_{j}^{-}}\right)\overline{N}_{j}\|\psi-\Omega\|_{\mathbb{B}},
\end{equation*}
therefore,
\begin{equation*}
\|\Gamma\psi-\Gamma\Omega\| \leq \kappa\|\psi-\Omega\|_{\mathbb{B}},
\,\ \text{where}\,\ \kappa<1.
\end{equation*}
According to the well-known contraction principle there exists a
unique fixed point $ h^{\ast }\left( \cdot \right) $ such that
$\Gamma h^{\ast }=h^{\ast }$. So, $h^{\ast }$ is a weighted
pseudo-almost automorphic on time-space scales solution of the model
(\ref{eq1}) in $\mathbb{F}=\{\psi\in\mathbb{B}: \,\
\|\psi\|_{\mathbb{B}}\leq r\}$. This completes the proof.
\begin{remark}
To the best of our knowledge, there have been no results focused on
the automorphic solutions, pseudo-almost automorphic ones and
Stepanov-like weighted pseudo-almost automorphic solutions on
time-space scales for high-order BAM neural networks with time
varying coefficients, mixed delays and leakage until now. Hence, the
obtained results are essentially new and the investigation methods
used in this paper can also be applied to study the Stepanov-like
weighted pseudo-almost automorphic solutions on time-space scales
for some other models of dynamical neural networks, such as
Cohen-Grossberg neural networks.
\end{remark}
\begin{theorem}\label{th1}
Let $(H_{1})-(H_{4})$ hold. The delta derivative of the only
Stepanov-like weighted pseudo-almost automorphic on time-space
scales solution of system (\ref{eq1}) is also Stepanov-like weighted
pseudo-almost automorphic on time-space scales (i.e. the unique
solution of (\ref{eq1}) is delta differentiable Stepanov-like
weighted pseudo-almost automorphic on time-space scales).
\end{theorem}
Proof. From system (\ref{eq1}), the expression of delta derivative
of the only Stepanov-like weighted pseudo-almost automorphic on
time-space scales solution is
\begin{eqnarray*}
x^{\Delta}_{i}\left( t\right)
&=&-\alpha_{i}(t)x_{i}\left(t-\eta_{i}(t)\right)
+\sum\limits_{j=1}^{m}D_{ij}\left( t\right) f_{j}\left(y_{j}\left( t\right) \right)\\
&+&\sum\limits_{j=1}^{m}D_{ij}^{\tau}\left( t\right)
f_{j}\left(y_{j}\left( t -\tau_{ij}(t)\right) \right)
+\sum\limits_{j=1}^{m}\overline{D}_{ij}\left( t\right)\int_{t-\sigma_{ij}(t)} ^{t}f_{j}\left(y_{j}\left( s\right) \right)\Delta s\\
&+&\sum\limits_{j=1}^{m}\sum\limits_{k=1}^{m}T_{ijk}(t)f_{k}(y_{k}(t-\chi_{k}(t)))f_{j}(y_{j}(t-\chi_{j}(t)))\\
&+&\sum\limits_{j=1}^{m}\widetilde{D}_{ij}\left(
t\right)\int_{t-\xi_{ij}(t)} ^{t}f_{j}\left(y_{j}^{\Delta}\left(
s\right) \right)\Delta s+I_{i}\left( t\right),
\end{eqnarray*}
and
\begin{eqnarray*}
y^{\Delta}_{j}\left( t\right)
&=&-c_{j}(t)y_{j}\left(t-\eta_{j}(t)\right)
+\sum\limits_{i=1}^{n}E_{ij}\left( t\right) f_{j}\left(x_{j}\left( t\right) \right)\\
&+&\sum\limits_{i=1}^{n}E_{ij}^{\tau}\left( t\right)
f_{j}\left(x_{j}\left( t -\tau_{ij}(t)\right) \right)
+\sum\limits_{i=1}^{n}\overline{E}_{ij}\left( t\right)\int_{t-\sigma_{ij}(t)}^{t}f_{j}\left(x_{j}\left( s\right) \right)\Delta s\\
&+&\sum\limits_{i=1}^{n}\sum\limits_{k=1}^{n}\overline{T}_{ijk}(t)f_{k}(x_{k}(t-\chi_{k}(t)))f_{j}(x_{j}(t-\chi_{j}(t)))\\
&+&\sum\limits_{i=1}^{n}\widetilde{D}_{ij}\left(
t\right)\int_{t-\xi_{ij}(t)}^{t}f_{j}\left(x_{j}^{\Delta}\left(
s\right) \right)\Delta s+J_{j}\left( t\right),\,\ t\in\mathbb{T}.
\end{eqnarray*}
Since all coefficients of the system (\ref{eq1}) are Stepanov-like
weighted pseudo-almost automorphic on time-space scales, derivative
of solution of system (\ref{eq1}) is Stepanov-like weighted
pseudo-almost automorphic on time-space scales.
\begin{remark}
In practice, time delays, leakage delay and parameter perturbations
are unavoidably encountered in the implementation of high-order BAM
neural networks, and they may destroy the stability of Stepanov-like
weighted pseudo-almost automorphic solution of HOBAMs, so it is
necessary and vital to study the dynamic behaviors of Stepanov-like
weighted pseudo-almost automorphic on time-space scales solution of
HOBAMs with time delays, leakage term and parameter perturbations.
\end{remark}
\section{Global exponential stability and convergence of Stepanov-like weighted pseudo-almost automorphic on time-space scales solution}
\subsection{Global exponential stability}
\begin{definition}
Let $Z^{\ast }\left(t\right) =\left(x_{1}^{\ast }\left(t\right),
x_{2}^{\ast }\left( t\right), \cdots ,x_{n}^{\ast
}\left(t\right),y_{1}^{\ast}\left(t\right),y_{2}^{\ast}\left(t\right),\cdots,
y_{m}^{\ast}\left(t\right) \right)^{T}$ be the Stepanov-like
weighted pseudo-almost automorphic solution on time-space scales of
system (\ref{eq1}) with initial value
$\psi^{\ast}\left(t\right)=\left(\varphi_{1}^{\ast}\left(t\right),\varphi_{2}^{\ast}\left(t\right),\cdots,\varphi_{n}^{\ast}\left(t\right),\right.$\\$\left.\phi_{1}^{\ast}\left(t\right),\phi_{2}^{\ast}\left(t\right),\cdots,
\phi_{m}^{\ast}\left(t\right)\right)^{T}$.
$Z^{\ast}\left(\cdot\right)$ is said to be globally exponential
stable if there exist constants $\gamma>0$,
$\ominus_{\nu}\gamma\in\mathcal{R}_{+}$ and $M>1$ such that for
every solution
\begin{equation*}
Z\left(t\right) =\left(x_{1}\left(t\right),x_{2}^{\ast
}\left(t\right),
\cdots,x_{n}\left(t\right),y_{1}\left(t\right),y_{2}\left(t\right),\cdots,y_{n}\left(t\right)
\right)^{T}
\end{equation*}
of system (\ref{eq1}) with any initial value
\begin{equation*}
\psi\left(t\right)=\left(\varphi_{1}\left(t\right),\varphi_{2}\left(t\right),\cdots,\varphi_{n}\left(t\right),\phi_{1}\left(t\right),\phi_{2}\left(t\right),\cdots,
\phi_{m}\left(t\right)\right)^{T}, \,\ \forall
t\in(0,+\infty)_{\mathbb{T}}, t\geq t_{0},
\end{equation*}
\begin{eqnarray*}
\left\Vert Z\left(t\right)-Z^{\ast}\left(t\right)\right\Vert_{0}
&=&\max \left\{\left\Vert
x\left(t\right)-x^{\ast}\left(t\right)\right\Vert_{\infty},\Vert
x^{\Delta}\left(t\right)
-x^{\ast\Delta}\left(t\right)\Vert_{\infty},\right.\\
&&\left.\left\Vert y\left(t\right)-y^{\ast
}\left(t\right)\right\Vert_{\infty},\Vert y^{\Delta}\left(t\right)
-y^{\ast\Delta }\left(t\right) \Vert_{\infty}\right\}\\
&\leq & M e_{\ominus\gamma}(t,t_{0})\|\psi\|_{1} \\
&=&M
e_{\ominus\gamma}(t,t_{0})\sup\limits_{t\in[-\theta,0]_{\mathbb{T}}}\max
\left\{\left\Vert\varphi\left(t\right)-\varphi^{\ast}\left(t\right)\right\Vert_{\infty},\left\Vert\varphi^{\Delta}\left(
t\right)-\varphi^{\ast\Delta}\left(t\right)
\right\Vert_{\infty},\right.\\
&&\left.\left\Vert\phi\left(t\right)-\phi^{\ast}\left(t\right)\right\Vert_{\infty},\left\Vert\phi^{\Delta}\left(
t\right)-\phi^{\ast\Delta}\left(t\right)
\right\Vert_{\infty}\right\},
\end{eqnarray*}
where $t_{0}=\max\{[-\theta,0]_{\mathbb{T}}\}$.
\end{definition}
\begin{theorem}\label{th2}
Let $(H_{1})-(H_{4})$ hold. The unique Stepanov-like weighted
pseudo-almost automorphic on time-space scales solution of system
(\ref{eq1}) is globally exponentially stable.
\end{theorem}
Proof. From Theorem \ref{th1} the system (\ref{eq1}) has one and
only one weighted pseudo-almost automorphic on time-space scales
solution on time scales
\begin{equation*}
Z^{*}(t)=(x^{*}_{1}(t),...x^{*}_{n}(t),y^{*}_{1}(t),...,y^{*}_{m}(t))^{T}\in\mathbb{R}^{n\times
m},
\end{equation*}
with the initial condition
\begin{equation*}
\psi^{*}(t)=(\varphi^{*}_{1}(t),...,\varphi^{*}_{n}(t),\phi^{*}_{1}(t),...,\phi^{*}_{m}(t))^{T}.
\end{equation*}
Let $Z(t)=(x_{1}(t),...,x_{n}(t),y_{1}(t),...,y_{m}(t))$ one arbitrary solution of (\ref{eq1}) with initial condition $\psi(t)=(\varphi_{1}(t),...,\varphi_{n}(t),\phi_{1}(t),...,\phi_{m}(t))^{T}$.\\

From system (\ref{eq1}), for $t\in\mathbb{T}$, we obtain:
\begin{equation}\label{eq4}
\left\{
\begin{array}{ccc}
u^{\Delta}_{i}\left( t\right) = -\alpha_{i}(t)u_{i}\left(t\right)
+\alpha_{i}(t)\int_{t-\eta_{i}(t)}^{t}u^{\Delta}_{i}(s)\Delta s
+\sum\limits_{j=1}^{m}D_{ij}\left( t\right)p_{j}\left(v_{j}(t)\right)\\
+\sum\limits_{j=1}^{m}D^{\tau}_{ij}\left(
t\right)p_{j}(v_{j}\left(\left(t-\tau_{ij}(t)\right)\right)
+\sum\limits_{j=1}^{m}\overline{D}_{ij}\left( t\right)\int_{t-\sigma_{ij}(t)} ^{t}p_{j}\left(v_{j}\left( s\right) \right)\Delta s\\
+\sum\limits_{j=1}^{m}\widetilde{D}_{ij}\left(
t\right)\int_{t-\xi_{ij}(t)}^{t}h_{j}\left(v_{j}^{\Delta}\left(
s\right) \right)\Delta s
+\sum\limits_{j=1}^{m}\sum\limits_{k=1}^{m}T_{ijk}(t)q_{j,k}(v_{j}(t-\chi_{j}(t)),v_{k}(t-\chi_{j}(t))),\\
v_{i}^{\Delta}(t)=-c_{i}(t)v_{i}(t)+c_{i}(t)\int_{t-\varsigma_{i}(t)}^{t}v_{i}^{\Delta}(s)\Delta
s+\sum\limits_{i=1}^{n}E_{ij}\left( t\right)p_{j}\left(u_{j}(t)\right)\\
+\sum\limits_{i=1}^{n}E^{\tau}_{ij}\left(
t\right)p_{j}(u_{j}\left(\left(t-\tau_{ij}(t)\right)\right)
+\sum\limits_{i=1}^{n}\overline{E}_{ij}\left( t\right)\int_{t-\sigma_{ij}(t)}^{t}p_{j}\left(u_{j}\left( s\right) \right)\Delta s\\
+\sum\limits_{i=1}^{n}\widetilde{E}_{ij}\left(
t\right)\int_{t-\xi_{ij}(t)}^{t}h_{j}\left(u_{j}^{\Delta}\left(
s\right)\right)\Delta s
+\sum\limits_{i=1}^{n}\sum\limits_{k=1}^{n}\overline{T}_{ijk}(t)q_{j,k}(u_{j}(t-\chi_{j}(t)),u_{k}(t-\chi_{j}(t))),\,\
\end{array}
\right.
\end{equation}
where
\begin{multline*} u_{i}(t)=x_{i}(t)-x^{*}_{i}(t),\,\ v_{i}(t)=y_{i}(t)-y_{i}^{*}(t),\,\
p_{j}(u_{j}(t))=f_{j}(x_{j}(t))-f_{j}(x_{j}^{*}(t)),\\ h_{j}(u_{j}^{\Delta}(t))=f_{j}(x_{j}^{\Delta}(t))-f_{j}({x_{j}^{*}}^{\Delta}(t)),\\
p_{j}(v_{j}(t))=f_{j}(y_{j}(t))-f_{j}(y_{j}^{*}(t)),\\ h_{j}(v_{j}^{\Delta}(t))=f_{j}(y_{j}^{\Delta}(t))-f_{j}({y_{j}^{*}}^{\Delta}(t)),\\
q_{j,k}(u_{j}(t),u_{k}(t))=f_{k}(x_{k}(t))f_{j}(x_{j}(t))-f_{k}(x^{*}_{k}(t))f_{j}(x^{*}_{j}(t))\\
q_{j,k}(v_{j}(t),v_{k}(t))=f_{k}(y_{k}(t))f_{j}(y_{j}(t))-f_{k}(y^{*}_{k}(t))f_{j}(y^{*}_{j}(t)).\\
\end{multline*}
For $i=1,...,n$ and $j=1,...,m$, the initial condition of
(\ref{eq4}) is
\begin{equation*}
u_{i}(s)=\varphi_{i}(s)-\varphi_{i}^{*}(s), \,\
v_{j}(s)=\phi_{j}(s)-\phi_{j}^{*}(s), \,\ s\in
[-\theta,0]_{\mathbb{T}}.
\end{equation*}

Multiplying the first equation in system (\ref{eq4}) by
$\hat{e}_{-\alpha_{i}}(t_{0},\sigma(s))$ and the second equation by
$\hat{e}_{-c_{j}}(t_{0},\sigma(s))$, and integrating over
$[t_{0},t]_{\mathbb{T}}$, where $t_{0}\in[-\theta,0]_{\mathbb{T}}$,
we obtain
\begin{equation}\label{eq5}
\left\{
\begin{array}{ccc}
u_{i}\left( t\right)=
u_{i}\left(t_{0}\right)\hat{e}_{-\alpha_{i}}(t,t_{0})+\int_{t_{0}}^{t}\hat{e}_{-\alpha_{i}}(t,\sigma(s))
\left(\alpha_{i}(s)\int_{s-\eta_{i}(s)}^{t}u^{\Delta}_{i}(u)\Delta u\right.\\
+\sum\limits_{j=1}^{m}D_{ij}\left(s\right)p_{j}\left(v_{j}(s)\right)
+\sum\limits_{j=1}^{m}D^{\tau}_{ij}\left(s\right)p_{j}(v_{j}\left(\left(s-\tau_{ij}(s)\right)\right)\\
+\sum\limits_{j=1}^{m}\overline{D}_{ij}\left(
s\right)\int_{s-\sigma_{ij}(s)}^{s}p_{j}\left(v_{j}\left( u\right)
\right)\Delta
u+\sum\limits_{j=1}^{m}\sum\limits_{k=1}^{m}T_{ijk}(s)q_{j,k}(v_{j}(s-\chi_{j}(s)),v_{k}(s-\chi_{k}(s)))
\\ \left.+\sum\limits_{j=1}^{m}\widetilde{D}_{ij}\left(
s\right)\int_{s-\xi_{ij}(s)}^{t}h_{j}\left(v_{j}^{\Delta}\left(
u\right) \right)\Delta u
\right)\Delta s,\\
v_{j}(t)=v_{j}(t_{0})\hat{e}_{-c_{j}}(t,t_{0})
+\int_{t_{0}}^{t}\hat{e}_{-c_{j}}(t,\sigma(s))\left(c_{j}(s)\int_{s-\varsigma_{j}(s)}^{s}v_{j}^{\Delta}(u)\Delta
u\right.\\
+\sum\limits_{i=1}^{n}E_{ij}\left(s\right)p_{i}\left(u_{i}(s)\right)
+\sum\limits_{i=1}^{n}E^{\tau}_{ij}\left(s\right)p_{i}(u_{i}\left(\left(s-\tau_{ij}(s)\right)\right)\\
+\sum\limits_{i=1}^{n}\overline{E}_{ij}\left(
s\right)\int_{s-\sigma_{ij}(s)}^{s}p_{i}\left(u_{i}\left( u\right)
\right)\Delta u
+\sum\limits_{i=1}^{n}\sum\limits_{k=1}^{n}\overline{T}_{ijk}(s)q_{j,k}(u_{i}(s-\chi_{i}(s)),u_{k}(s-\chi_{k}(s)))\\
\left.\left.+\sum\limits_{i=1}^{n}\widetilde{E}_{ij}\left(
s\right)\int_{s-\xi_{ij}(s)}^{t}h_{i}\left(u_{i}^{\Delta}\left(
u\right) \right)\Delta u \right)\Delta s\right)\Delta s,\,\
t\in\mathbb{T},
\end{array}
\right.
\end{equation}
Now, we define $G_{i}$, $\overline{G}_{j}$, $H_{i}$ and $\overline{H}_{j}$ as follows:\\
\begin{multline*}
G_{i}(w)=\alpha_{i}^{-}-w-\left(\exp\left(w
\sup\limits_{s\in\mathbb{T}}\nu(s)\right)\left(\alpha_{i}^{+}\eta_{i}^{+}
\exp\left(w\eta_{i}^{+}\right)+\sum\limits_{j=1}^{m}D_{ij}^{+}L_{j}\right.\right.\\
+\sum\limits_{j=1}^{m}(D_{ij}^{\tau})^{+}L_{j}\exp\left(w\tau_{ij}^{+}\right)
+\sum\limits_{j=1}^{m}\overline{D}_{ij}^{+}L_{j}\sigma_{ij}^{+}\exp\left(w\sigma_{ij}^{+}\right)
+\sum\limits_{j=1}^{m}\widetilde{D}_{ij}^{+}L_{j}\xi_{ij}^{+}\exp\left(w\xi_{ij}^{+}\right)\\
\left.+\sum\limits_{j=1}^{m}\sum\limits_{k=1}^{m}T_{ijk}^{+}(L_{k}r+|f_{k}(0)|)(L_{j}r+|f_{j}(0)|)\right),\phantom{++++++}\\
\overline{G}_{j}(w)=c_{j}^{-}-w-\left(\exp\left(w
\sup\limits_{s\in\mathbb{T}}\nu(s)\right)\left(c_{j}^{+}\varsigma_{j}^{+}
\exp\left(w\varsigma_{j}^{+}\right)+\sum\limits_{i=1}^{n}E_{ij}^{+}L_{i}\right.\right.\\
+\sum\limits_{i=1}^{n}(E_{ij}^{\tau})^{+}L_{i}\exp\left(w\tau_{ij}^{+}\right)
+\sum\limits_{i=1}^{n}\overline{E}_{ij}^{+}L_{i}\sigma_{ij}^{+}\exp\left(w\sigma_{ij}^{+}\right)
+\sum\limits_{i=1}^{n}\widetilde{E}_{ij}^{+}L_{i}\xi_{ij}^{+}\exp\left(w\xi_{ij}^{+}\right)\\
\left.+\sum\limits_{i=1}^{n}\sum\limits_{k=1}^{n}\overline{T}_{ijk}^{+}(L_{k}r+|f_{k}(0)|)(L_{i}r+|f_{i}(0)|)\right),\phantom{++++++}\\
 H_{i}(w)=\alpha_{i}^{-}-w-\left(\alpha_{i}^{+}\exp\left(w
\sup\limits_{s\in\mathbb{T}}\nu(s)+\alpha_{i}^{-}-\beta\right)\left(\alpha_{i}^{+}\eta_{i}^{+}
\exp\left(w\eta_{i}^{+}\right)\right.\right.\\
+\sum\limits_{j=1}^{m}D_{ij}^{+}L_{j}+\sum\limits_{j=1}^{m}(D_{ij}^{\tau})^{+}L_{j}\exp\left(w\tau_{ij}^{+}\right)
+\sum\limits_{j=1}^{m}\overline{D}_{ij}^{+}L_{j}\sigma_{ij}^{+}\exp\left(w\sigma_{ij}^{+}\right)
+\sum\limits_{j=1}^{m}\widetilde{D}_{ij}^{+}L_{j}\xi_{ij}^{+}\exp\left(w\xi_{ij}^{+}\right)
\end{multline*}
\begin{multline*}
\left.+\sum\limits_{j=1}^{m}\sum\limits_{k=1}^{m}T_{ijk}^{+}(L_{k}r+|f_{k}(0)|)(L_{j}r+|f_{j}(0)|)\right),\phantom{++++++}\\
\overline{H}_{j}(w)=c_{j}^{-}-w-\left(c_{j}^{+}\exp\left(w
\sup\limits_{s\in\mathbb{T}}\nu(s)+c_{j}^{-}-\beta\right)\left(c_{j}^{+}\varsigma_{j}^{+}
\exp\left(w\varsigma_{j}^{+}\right)\right.\right.\\
+\sum\limits_{i=1}^{n}E_{ij}^{+}L_{i}+\sum\limits_{i=1}^{n}(E_{ij}^{\tau})^{+}L_{i}\exp\left(w\tau_{ij}^{+}\right)
+\sum\limits_{i=1}^{n}\overline{E}_{ij}^{+}L_{i}\sigma_{ij}^{+}\exp\left(w\sigma_{ij}^{+}\right)
+\sum\limits_{i=1}^{n}\widetilde{E}_{ij}^{+}L_{i}\xi_{ij}^{+}\exp\left(w\xi_{ij}^{+}\right)\\
\left.+\sum\limits_{i=1}^{n}\sum\limits_{k=1}^{n}T_{ijk}^{+}(L_{k}r+|f_{k}(0)|)(L_{i}r+|f_{i}(0)|)\right),\phantom{++++++}\\
\end{multline*}
where $i=1,...,n$, $j=1,...,m$, $w\in(0,+\infty)$.\\

From $\left( H_{3} \right)$, we have\\
\begin{multline*}G_{i}(0)=\alpha_{i}^{-}-\left(\alpha_{i}^{+}\eta_{i}^{+}
+\sum\limits_{j=1}^{m}D_{ij}^{+}L_{j}+\sum\limits_{j=1}^{m}(D_{ij}^{\tau})^{+}L_{j}
+\sum\limits_{j=1}^{m}\overline{D}_{ij}^{+}L_{j}\sigma_{ij}^{+}\right.\\
\left.+\sum\limits_{j=1}^{m}\widetilde{D}_{ij}^{+}L_{j}\xi_{ij}^{+}
+\sum\limits_{j=1}^{m}\sum\limits_{k=1}^{m}T_{ijk}^{+}(L_{k}r+|f_{k}(0)|)(L_{j}r+|f_{j}(0)|)\right),\\
\overline{G}_{j}(0)=c_{j}^{-}+\sum\limits_{i=1}^{n}E_{ij}^{+}L_{i}
+\sum\limits_{i=1}^{n}(E_{ij}^{\tau})^{+}L_{i}
+\sum\limits_{i=1}^{n}\overline{E}_{ij}^{+}L_{i}\sigma_{ij}^{+}\\
\left.+\sum\limits_{i=1}^{n}\widetilde{E}_{ij}^{+}L_{i}\xi_{ij}^{+}
+\sum\limits_{i=1}^{n}\sum\limits_{k=1}^{n}\overline{T}_{ijk}^{+}(L_{k}r+|f_{k}(0)|)(L_{i}r+|f_{i}(0)|)\right),\\
H_{i}(0)=\alpha_{i}^{-}-\alpha_{i}^{+}\exp\left(\alpha_{i}^{-}-\beta\right)\left(\alpha_{i}^{+}\eta_{i}^{+}
+\sum\limits_{j=1}^{m}D_{ij}^{+}L_{j}+\sum\limits_{j=1}^{m}(D_{ij}^{\tau})^{+}L_{j}
+\sum\limits_{j=1}^{m}\overline{D}_{ij}^{+}L_{j}\sigma_{ij}^{+}\right.\\
\left.+\sum\limits_{j=1}^{m}\widetilde{D}_{ij}^{+}L_{j}\xi_{ij}^{+}
+\sum\limits_{j=1}^{m}\sum\limits_{k=1}^{m}T_{ijk}^{+}(L_{k}r+|f_{k}(0)|)(L_{j}r+|f_{j}(0)|)\right),
\end{multline*}
\begin{multline*}
\overline{H}_{j}(0)=c_{j}^{-}-c_{j}^{+}\exp\left(c_{j}^{-}-\beta\right)\left(c_{j}^{+}\varsigma_{j}^{+}
+\sum\limits_{i=1}^{n}E_{ij}^{+}L_{i}+\sum\limits_{i=1}^{n}(E_{ij}^{\tau})^{+}L_{i}
+\sum\limits_{i=1}^{n}\overline{E}_{ij}^{+}L_{i}\sigma_{ij}^{+}\right.\\
\left.+\sum\limits_{i=1}^{n}\widetilde{E}_{ij}^{+}L_{i}\xi_{ij}^{+}
+\sum\limits_{i=1}^{n}\sum\limits_{k=1}^{n}\overline{T}_{ijk}^{+}(L_{k}r+|f_{k}(0)|)(L_{i}r+|f_{i}(0)|)\right),
\end{multline*}
Since the functions $G_{i}(.)$, $\overline{G}_{j}(.)$, $H_{i}(.)$
and $\overline{H}_{j}(.)$ are continuous on $[0,+\infty)$ and
$G_{i}(w)$, $\overline{G}_{j}(w)$, $H_{i}(w)$,
$\overline{H}_{j}(w)\longrightarrow-\infty$ when
$w\longrightarrow+\infty$, it exist $\eta_{i}, \bar{\eta}_{j},
\epsilon_{i}, \bar{\epsilon}_{j}>0$ such as
\begin{equation*}
H_{i}(\eta_{i})=\overline{H}_{j}(\bar{\eta}_{j})=G_{i}(\epsilon_{i})=\overline{G}_{j}(\bar{\epsilon}_{j})=0
\end{equation*}
and
\begin{eqnarray*}
&&G_{i}(w)>0\,\ \text{for} \,\ w\in(0,\eta_{i}), \,\ \overline{G}_{j}(w)>0 \,\ \text{for} \,\ w\in(0,\bar{\eta}_{j}),\\
&&H_{i}(w)>0\,\ \text{for} \,\ w\in(0,\epsilon_{i}), \,\
\overline{H}_{j}(w)>0\,\ \text{for} \,\ w\in(0,\bar{\epsilon}_{j}).
\end{eqnarray*}

Let $a=\min\limits_{1\leq i\leq n}\left\{\eta_{i}, \bar{\eta}_{j},
\epsilon_{i}, \bar{\epsilon}_{j}\right\}$, we obtain
\begin{equation*}
H_{i}(a)\geq0,\,\ \overline{H}_{j}(a)\geq0, \,\ G_{i}(a)\geq0,\,\
\text{and} \,\ \overline{G}_{j}(a)\geq0,\,\ i=1,...,n, j=1,...,m.
\end{equation*}
So, we can choose the positive constant
$0<\gamma<\min\limits_{1\leq i \leq n, 1\leq j\leq m}\{a,\alpha_{i}^{-},c_{j}^{-}\},$\\
$\text{such that}\,\ H_{i}(\gamma)>0,\,\ \overline{H}_{j}(\gamma)>0,
\,\ G_{i}(\gamma)>0\,\ \text{and} \,\ \overline{G}_{j}(\gamma)>0,
\,\ i=1,...,n, j=1,...,m.$ which imply that, for $i=1,...,n$ and
$j=1,...,m$
\begin{multline*}
\frac{1}{\alpha_{i}^{-}-\gamma}\sum\limits_{j=1}^{m}\left(\exp(\gamma\sup\limits_{s\in\mathbb{T}}\nu(s))\left(\alpha_{i}^{+}\eta_{i}^{+}+\sum\limits_{j=1}^{m}
\left(D_{ij}^{+}+(D_{ij}^{\tau})^{+}+\overline{D}_{ij}^{+}\sigma_{ij}^{+}
+\widetilde{D}_{ij}^{+}\xi_{ij}^{+}\right)L_{j}\right.\right.\\
\left.\left.+\sum\limits_{j=1}^{m}\sum\limits_{k=1}^{n}T_{ijk}^{+}(L_{k}r+|f_{k}(0)|)(L_{j}r+|f_{j}(0)|)\right)\right)<1,\\
\frac{1}{c_{j}^{-}-\gamma}\sum\limits_{i=1}^{n}\left(\exp(\gamma\sup\limits_{s\in\mathbb{T}}\nu(s))\left(c_{j}^{+}\varsigma_{j}^{+}+\sum\limits_{i=1}^{n}
\left(E_{ij}^{+}+(E_{ij}^{\tau})^{+}+\overline{E}_{ij}^{+}\sigma_{ij}^{+}
+\widetilde{E}_{ij}^{+}\xi_{ij}^{+}\right)L_{i}\right.\right.\\
\left.\left.+\sum\limits_{i=1}^{n}\sum\limits_{k=1}^{n}\overline{T}_{ijk}^{+}(L_{k}r+|f_{k}(0)|)(L_{i}r+|f_{i}(0)|)\right)\right)<1,\\
\left(1+\frac{\alpha_{i}^{+}\exp(\gamma\sup\limits_{s\in\mathbb{T}}\nu(s))}{\alpha_{i}^{-}-\gamma}\right)
\sum\limits_{j=1}^{m}\left(\exp(\gamma\sup\limits_{s\in\mathbb{T}}\nu(s))\left(\alpha_{i}^{+}\eta_{i}^{+}
+\sum\limits_{j=1}^{m}\left(D_{ij}^{+}+(D_{ij}^{\tau})^{+}\right.\right.\right.\\
\left.\left.\left.+\overline{D}_{ij}^{+}\sigma_{ij}^{+}+\widetilde{D}_{ij}^{+}\xi_{ij}^{+}\right)L_{j}
+\sum\limits_{j=1}^{m}\sum\limits_{k=1}^{m}T_{ijk}^{+}(L_{k}r+|f_{k}(0)|)(L_{j}r+|f_{j}(0)|)\right)\right)<1,\\
\left(1+\frac{c_{j}^{+}\exp(\gamma\sup\limits_{s\in\mathbb{T}}\nu(s))}{c_{j}^{-}-\gamma}\right)
\sum\limits_{i=1}^{n}\left(\exp(\gamma\sup\limits_{s\in\mathbb{T}}\nu(s))\left(c_{j}^{+}\varsigma_{j}^{+}
+\sum\limits_{i=1}^{n}\left(E_{ij}^{+}+(E_{ij}^{\tau})^{+}\right.\right.\right.\\
\left.\left.\left.+\overline{E}_{ij}^{+}\sigma_{ij}^{+}+\widetilde{E}_{ij}^{+}\xi_{ij}^{+}\right)L_{i}
+\sum\limits_{i=1}^{n}\sum\limits_{k=1}^{n}\overline{T}_{ijk}^{+}(L_{k}r+|f_{k}(0)|)(L_{i}r+|f_{i}(0)|)\right)\right)<1,\\
\end{multline*}
Let
\begin{equation}\label{abcd5}
K=\max\limits_{1\leq i\leq n, 1\leq j\leq
m}\left\{\frac{\alpha_{i}^{-}}{K^{*}},
\frac{c_{j}^{-}}{P^{*}}\right\},
\end{equation}
where
\begin{multline*}
\phantom{++++++}K^{*}=\alpha_{i}^{+}\eta_{i}^{+}+\sum\limits_{j=1}^{m}\left(D_{ij}^{+}+(D_{ij}^{\tau})^{+}+\overline{D}_{ij}^{+}\sigma_{ij}^{+}
+\widetilde{D}_{ij}^{+}\xi_{ij}^{+}\right)L_{j}\\
+\sum\limits_{j=1}^{m}\sum\limits_{k=1}^{m}T_{ijk}^{+}(L_{k}r+|f_{k}(0)|)(L_{j}r+|f_{j}(0)|),\phantom{+++}
\end{multline*}
and
\begin{multline*}
\phantom{++++++}P^{*}=c_{j}^{+}\varsigma_{j}^{+}+\sum\limits_{i=1}^{n}\left(E_{ij}^{+}+(E_{ij}^{\tau})^{+}+\overline{E}_{ij}^{+}\sigma_{ij}^{+}
+\widetilde{E}_{ij}^{+}\xi_{ij}^{+}\right)L_{i}\\
+\sum\limits_{i=1}^{n}\sum\limits_{k=1}^{n}\overline{T}_{ijk}^{+}(L_{k}r+|f_{k}(0)|)(L_{i}r+|f_{i}(0)|).\phantom{+++}
\end{multline*}
By hypothesis ($H_{3}$), we have $K>1$, therefore,
\begin{equation}\label{abcd6}
\|Z(t)-Z^{*}(t)\|\leq K \hat{e}_{\ominus_{\nu}
\gamma}(t,t_{0})\|\psi\|_{0},\,\ \forall t\in[t_{0},0]_{\mathbb{T}},
\end{equation}
where $\ominus_{\nu}\gamma\in \mathcal{R}_{\nu}^{+}$. We claim that
\begin{equation}\label{abcd71}
\|Z(t)-Z^{*}(t)\|\leq K \hat{e}_{\ominus_{\nu}
\gamma}(t,t_{0})\|\psi\|_{0},\,\ \forall
t\in[t_{0},+\infty)_{\mathbb{T}}.
\end{equation}
To prove (\ref{abcd71}), we show that for any $\varpi>1$, the
following inequality holds:
\begin{equation}\label{abcd72}
\|Z(t)-Z^{*}(t)\|\leq \varpi K \hat{e}_{\ominus_{\nu}
\gamma}(t,t_{0})\|\psi\|_{0},\,\ \forall
t\in[t_{0},+\infty)_{\mathbb{T}}.
\end{equation}
If (\ref{abcd72}) is not true, then there must be some
$t_{1}\in(0,+\infty)_{\mathbb{T}}$, $d\geq1$ such that
\begin{equation}\label{abcd73}
\|Z(t_{1})-Z^{*}(t_{1})\|=d \varpi K \hat{e}_{\ominus_{\nu}
\gamma}(t_{1},t_{0})\|\psi\|_{0},
\end{equation}
and
\begin{equation}\label{abcd74}
\|Z(t)-Z^{*}(t)\|\leq d \varpi K \hat{e}_{\ominus_{\nu}
\gamma}(t,t_{0})\|\psi\|_{0}, \,\ t\in[t_{0},t_{1}]_{\mathbb{T}}.
\end{equation}
By (\ref{eq5}), (\ref{abcd73}), (\ref{abcd74}) and
$(H_{1})-(H_{3})$, we have for $i=1,...,n$
\begin{eqnarray*}
|u_{i}(t_{1})|&\leq&
\hat{e}_{-\alpha_{i}}(t_{1},t_{0})\|\psi\|_{0}+d \varpi K
\hat{e}_{\ominus_{\nu}\gamma}(t_{1},t_{0})\|\psi\|_{0}\int_{t_{0}}^{t_{1}}
\hat{e}_{-\alpha_{i}}(t_{1},\sigma(s))\hat{e}_{\gamma}(t_{1},\sigma(s))\nonumber\\
&\times&\left(\alpha_{i}^{+}\int_{s-\eta_{i}(s)}^{s}\hat{e}_{\gamma}(\sigma(u),u)\Delta
u+\sum\limits_{j=1}^{n}
D_{ij}^{+}L_{j}\hat{e}_{\gamma}(\sigma(s),s)\right.\nonumber\\
&+&\sum\limits_{j=1}^{n}\sum\limits_{k=1}^{n}T_{ijk}^{+}(L_{k}r+|f_{k}(0)|)(L_{j}r+|f_{j}(0)|)\hat{e}_{\gamma}(\sigma(s),s-\chi_{j}(s))\hat{e}_{\gamma}(\sigma(s),
s-\chi_{k}(s))\nonumber\\
&+&(D_{ij}^{\tau})^{+}L_{j}\hat{e}_{\gamma}(\sigma(s),s-\tau_{ij}(s))+\sum\limits_{j=1}^{n}
\overline{D}_{ij}^{+}L_{j}\int_{s-\sigma_{ij}(s)}^{s}\hat{e}_{\gamma}(\sigma(u),u)\Delta u\nonumber\\
&+&\left. \sum\limits_{j=1}^{n}\widetilde{D}_{ij}^{+}L_{j}\int_{s-\xi_{ij}(s)}^{s}\hat{e}_{\gamma}(\sigma(u),u)\Delta u\right)\Delta s\nonumber\\
&\leq& \hat{e}_{-\alpha_{i}}(t_{1},t_{0})\|\psi\|_{0}+d \varpi K
\hat{e}_{\ominus_{\nu}\gamma}(t_{1},t_{0})\|\psi\|_{0}\int_{t_{0}}^{t_{1}}
\hat{e}_{-\alpha_{i}}(t_{1},\sigma(s))\hat{e}_{\gamma}(t_{1},\sigma(s))\nonumber\\
&\times&\left(\alpha_{i}^{+}\hat{e}_{\gamma}(\sigma(s),s-\eta_{i}(s))+\sum\limits_{j=1}^{n}
D_{ij}^{+}L_{j}\hat{e}_{\gamma}(\sigma(s),s)\right.\nonumber\\
&+&\sum\limits_{j=1}^{n}\sum\limits_{k=1}^{n}T_{ijk}^{+}(L_{k}r+|f_{k}(0)|)(L_{j}r+|f_{j}(0)|)\hat{e}_{\gamma}(\sigma(s),s-\chi_{j}(s))\hat{e}_{\gamma}(\sigma(s)
,s-\chi_{k}(s))\nonumber\\
&+&\sum\limits_{j=1}^{n}(D_{ij}^{\tau})^{+}L_{j}\hat{e}_{\gamma}(\sigma(s),s-\tau_{ij}(s))+\sum\limits_{j=1}^{n}
\overline{D}_{ij}^{+}L_{j}\int_{s-\sigma_{ij}(s)}^{s}\hat{e}_{\gamma}(\sigma(u),u)\Delta
u \nonumber\\&+&\left. \sum\limits_{j=1}^{n}
\widetilde{D}_{ij}^{+}L_{j}\int_{s-\xi_{ij}(s)}^{s}\hat{e}_{\gamma}(\sigma(u),u)\Delta u\right)\Delta s\nonumber\\
&\leq& \hat{e}_{-\alpha_{i}}(t_{1},t_{0})\|\psi\|_{0}+d \varpi K
\hat{e}_{\ominus_{\nu}\gamma}(t_{1},t_{0})\|\psi\|_{0}\int_{t_{0}}^{t_{1}}
\hat{e}_{-\alpha_{i}}(t_{1},\sigma(s))\hat{e}_{\gamma}(t_{1},\sigma(s))\nonumber\\
&\times&\left(\alpha_{i}^{+}\eta_{i}^{+}\exp\left[\gamma\left(\eta_{i}^{+}+\sup\limits_{s\in\mathbb{T}}\nu(s)\right)\right]
+\sum\limits_{j=1}^{n}D_{ij}^{+}L_{j}\exp\left(\gamma\sup\limits_{s\in\mathbb{T}}\nu(s)\right)\right.\nonumber\\
&+&(D_{ij}^{\tau})^{+}L_{j}\exp\left[\gamma\left(\tau_{ij}^{+}+\sup\limits_{s\in\mathbb{T}}\nu(s)\right)\right]
+\sum\limits_{j=1}^{n}\overline{D}_{ij}^{+}L_{j}\exp\left[\gamma\left(\sigma_{ij}^{+}+\sup\limits_{s\in\mathbb{T}}\nu(s)\right)\right]\\
&+&\sum\limits_{j=1}^{n}\widetilde{D}_{ij}^{+}L_{j}\exp\left[\gamma\left(\xi_{ij}^{+}+\sup\limits_{s\in\mathbb{T}}\nu(s)\right)\right]
+\sum\limits_{j=1}^{n}\sum\limits_{k=1}^{n}T_{ijk}^{+}(L_{k}r+|f_{k}(0)|)(L_{j}r+|f_{j}(0)|)\nonumber\\
&\times&\left.\exp\left[\gamma\left(\chi_{j}^{+}+\sup\limits_{s\in\mathbb{T}}\nu(s)\right)\right]\exp\left[\gamma\left(\chi_{k}^{+}
+\sup\limits_{s\in\mathbb{T}}\nu(s)\right)\right] \right)\Delta s
\end{eqnarray*}
\begin{eqnarray}\label{from}
&\leq& d\varpi K \hat{e}_{\ominus_{\nu}
\gamma}(t_{1},t_{0})\|\psi\|_{0}\left\{\frac{1}{K}\hat{e}_{-\alpha_{i}\oplus_{\nu}\gamma}(t_{1},t_{0})+
\left[\exp\left(\gamma\sup\limits_{s\in\mathbb{T}}\nu(s)\right)\right.\right.\nonumber\\
&\times& \left(\alpha_{i}^{+}\eta_{i}^{+}\exp(\gamma\eta_{i}^{+})+
\sum\limits_{j=1}^{n}D_{ij}^{+}L_{j}+\sum\limits_{j=1}^{n}(D_{ij}^{\tau})^{+}L_{j}\exp(\gamma\tau_{ij}^{+})
+\sum\limits_{j=1}^{n}\overline{D}_{ij}^{+}\sigma_{ij}^{+}L_{j}\exp(\gamma\sigma_{ij}^{+})\right.\nonumber\\
&+&\sum\limits_{j=1}^{n}\widetilde{D}_{ij}^{+}L_{j}\xi_{ij}^{+}\exp(\gamma\xi_{ij}^{+})
+\sum\limits_{j=1}^{n}\sum\limits_{k=1}^{n}T_{ijk}^{+}(L_{k}r+|f_{k}(0)|)(L_{j}r+|f_{j}(0)|)\nonumber\\
&\times&\left.\left.\left.\exp(\gamma \chi_{j}^{+})\exp(\gamma
\chi_{k}^{+})\right)\right]\frac{1-\hat{e}_{-\alpha_{i}\oplus_{\nu}\gamma}(t_{1},t_{0})}{\alpha_{i}^{-}-\gamma}\right\}\nonumber\\
&\leq& d\varpi K \hat{e}_{\ominus_{\nu}
\gamma}(t_{1},t_{0})\|\psi\|_{0}\left\{\frac{1}{K}-\frac{1}{\alpha_{i}^{-}-\gamma}\left(
\exp\left(\gamma\sup\limits_{s\in\mathbb{T}}\nu(s)\right) \left(\alpha_{i}^{+}\eta_{i}^{+}\exp(\gamma\eta_{i}^{+})\right.\right.\right.\nonumber\\
&+&\sum\limits_{j=1}^{n}D_{ij}^{+}L_{j}+\sum\limits_{j=1}^{n}(D_{ij}^{\tau})^{+}L_{j}\exp(\gamma\tau_{ij}^{+})
+\sum\limits_{j=1}^{n}\overline{D}_{ij}^{+}\sigma_{ij}^{+}L_{j}\exp(\gamma\sigma_{ij}^{+})\nonumber\\
&+&\sum\limits_{j=1}^{n}\widetilde{D}_{ij}^{+}L_{j}\xi_{ij}^{+}\exp(\gamma\xi_{ij}^{+})
+\sum\limits_{j=1}^{n}\sum\limits_{k=1}^{n}T_{ijk}^{+}(L_{k}r+|f_{k}(0)|)(L_{j}r+|f_{j}(0)|)\nonumber\\
&\times&\left.\left.\left.\exp(\gamma \chi_{j}^{+})\exp(\gamma
\chi_{k}^{+})\right)\right) +
B_{i}^{+}\right]\hat{e}_{-\alpha_{i}\oplus_{\nu}\gamma}(t_{1},t_{0})
+\frac{1}{\alpha_{i}^{-}-\gamma}\left(\exp\left(\gamma\sup\limits_{s\in\mathbb{T}}\nu(s)
\left(\alpha_{i}^{+}\eta_{i}^{+}\exp(\gamma\eta_{i}^{+})\right.\right.\right.\nonumber\\
&+&\sum\limits_{j=1}^{n}D_{ij}^{+}L_{j}+\sum\limits_{j=1}^{n}(D_{ij}^{\tau})^{+}L_{j}\exp(\gamma\tau_{ij}^{+})
+\sum\limits_{j=1}^{n}\overline{D}_{ij}^{+}\sigma_{ij}^{+}L_{j}\exp(\gamma\sigma_{ij}^{+})\nonumber\\
&+&\sum\limits_{j=1}^{n}\widetilde{D}_{ij}^{+}L_{j}\xi_{ij}^{+}\exp(\gamma\xi_{ij}^{+})
+\sum\limits_{j=1}^{n}\sum\limits_{k=1}^{n}T_{ijk}^{+}(L_{k}r+|f_{k}(0)|)(L_{j}r+|f_{j}(0)|)\nonumber\\
&\times&\left.\left.\left.\left.\exp(\gamma \chi_{j}^{+})\exp(\gamma \chi_{k}^{+})\right)\right)\right)\right\}\nonumber\\
&\leq& d\varpi K \hat{e}_{\ominus_{\nu}
\gamma}(t_{1},t_{0})\|\psi\|_{0}.
\end{eqnarray}
In addition, we have
\begin{eqnarray*}
|v_{j}(t_{1})| &\leq& \hat{e}_{-c_{j}}(t_{1},t_{0})\|\psi\|_{0}+d
\varpi K
\hat{e}_{\ominus_{\nu}\gamma}(t_{1},t_{0})\|\psi\|_{0}\int_{t_{0}}^{t_{1}}
\hat{e}_{-c_{j}}(t_{1},\sigma(s))\hat{e}_{\gamma}(t_{1},\sigma(s))\nonumber\\
&\times&\left(c_{j}^{+}\varsigma_{j}^{+}\exp\left[\gamma\left(\varsigma_{j}^{+}+\sup\limits_{s\in\mathbb{T}}\nu(s)\right)\right]
+\sum\limits_{i=1}^{n}E_{ij}^{+}L_{i}\exp\left(\gamma\sup\limits_{s\in\mathbb{T}}\nu(s)\right)\right.\nonumber\\
&+&(E_{ij}^{\tau})^{+}L_{i}\exp\left[\gamma\left(\tau_{ij}^{+}+\sup\limits_{s\in\mathbb{T}}\nu(s)\right)\right]
+\sum\limits_{i=1}^{n}\overline{E}_{ij}^{+}L_{i}\exp\left[\gamma\left(\sigma_{ij}^{+}+\sup\limits_{s\in\mathbb{T}}\nu(s)\right)\right]
\\&+&\sum\limits_{i=1}^{n}\widetilde{E}_{ij}^{+}L_{i}\exp\left[\gamma\left(\xi_{ij}^{+}+\sup\limits_{s\in\mathbb{T}}\nu(s)\right)\right]
+\sum\limits_{i=1}^{n}\sum\limits_{k=1}^{n}\overline{T}_{ijk}^{+}(L_{k}r+|f_{k}(0)|)(L_{i}r+|f_{i}(0)|)\nonumber\\
&\times&\left.\exp\left[\gamma\left(\chi_{i}^{+}+\sup\limits_{s\in\mathbb{T}}\nu(s)\right)\right]\exp\left[\gamma\left(\chi_{k}^{+}
+\sup\limits_{s\in\mathbb{T}}\nu(s)\right)\right]
\right)\Delta s\nonumber\\
\end{eqnarray*}
\begin{eqnarray*}
&\leq& d\varpi K \hat{e}_{\ominus_{\nu}
\gamma}(t_{1},t_{0})\|\psi\|_{0}\left\{\frac{1}{K}\hat{e}_{-c_{j}\oplus_{\nu}\gamma}(t_{1},t_{0})+
\left[\exp\left(\gamma\sup\limits_{s\in\mathbb{T}}\nu(s)\right)\right.\right.\nonumber\\
&\times&
\left(c_{j}^{+}\varsigma_{j}^{+}\exp(\gamma\varsigma_{j}^{+})+
\sum\limits_{i=1}^{n}E_{ij}^{+}L_{i}+\sum\limits_{i=1}^{n}(E_{ij}^{\tau})^{+}L_{i}\exp(\gamma\tau_{ij}^{+})
+\sum\limits_{i=1}^{n}\overline{E}_{ij}^{+}\sigma_{ij}^{+}L_{i}\exp(\gamma\sigma_{ij}^{+})\right.\nonumber\\
&+&\sum\limits_{i=1}^{n}\widetilde{E}_{ij}^{+}L_{i}\xi_{ij}^{+}\exp(\gamma\xi_{ij}^{+})
+\sum\limits_{i=1}^{n}\sum\limits_{k=1}^{n}\overline{T}_{ijk}^{+}(L_{k}r+|f_{k}(0)|)(L_{i}r+|f_{i}(0)|)\nonumber\\
&\times&\left.\left.\left.\exp(\gamma \chi_{i}^{+})\exp(\gamma
\chi_{k}^{+})\right)\right]\frac{1-\hat{e}_{-c_{j}\oplus_{\nu}\gamma}(t_{1},t_{0})}{c_{j}^{-}-\gamma}\right\}\nonumber\\
&\leq& d\varpi K \hat{e}_{\ominus_{\nu}
\gamma}(t_{1},t_{0})\|\psi\|_{0}\left\{\frac{1}{K}-\frac{1}{c_{j}^{-}-\gamma}\left(
\exp\left(\gamma\sup\limits_{s\in\mathbb{T}}\nu(s)\right) \left(c_{j}^{+}\varsigma_{j}^{+}\exp(\gamma\varsigma_{j}^{+})\right.\right.\right.\nonumber\\
&+&\sum\limits_{i=1}^{n}E_{ij}^{+}L_{i}+\sum\limits_{i=1}^{n}(E_{ij}^{\tau})^{+}L_{i}\exp(\gamma\tau_{ij}^{+})
+\sum\limits_{i=1}^{n}\overline{E}_{ij}^{+}\sigma_{ij}^{+}L_{i}\exp(\gamma\sigma_{ij}^{+})\nonumber\\
&+&\sum\limits_{i=1}^{n}\widetilde{E}_{ij}^{+}L_{i}\xi_{ij}^{+}\exp(\gamma\xi_{ij}^{+})
+\sum\limits_{i=1}^{n}\sum\limits_{k=1}^{n}\overline{T}_{ijk}^{+}(L_{k}r+|f_{k}(0)|)(L_{i}r+|f_{i}(0)|)\nonumber\\
&\times&\left.\left.\left.\exp(\gamma \chi_{i}^{+})\exp(\gamma
\chi_{k}^{+})\right)\right)\right]\hat{e}_{-c_{j}\oplus_{\nu}\gamma}(t_{1},t_{0})
+\frac{1}{c_{j}^{-}-\gamma}\left(\exp\left(\gamma\sup\limits_{s\in\mathbb{T}}\nu(s)
\left(c_{j}^{+}\varsigma_{j}^{+}\exp(\gamma\varsigma_{j}^{+})\right.\right.\right.\nonumber\\
&+&\sum\limits_{i=1}^{n}E_{ij}^{+}L_{i}+\sum\limits_{i=1}^{n}(E_{ij}^{\tau})^{+}L_{i}\exp(\gamma\tau_{ij}^{+})
+\sum\limits_{i=1}^{n}\overline{E}_{ij}^{+}\sigma_{ij}^{+}L_{i}\exp(\gamma\sigma_{ij}^{+})\nonumber\\
&+&\sum\limits_{i=1}^{n}\widetilde{E}_{ij}^{+}L_{i}\xi_{ij}^{+}\exp(\gamma\xi_{ij}^{+})
+\sum\limits_{i=1}^{n}\sum\limits_{k=1}^{n}\overline{T}_{ijk}^{+}(L_{k}r+|f_{k}(0)|)(L_{i}r+|f_{i}(0)|)\nonumber\\
&\times&\left.\left.\left.\left.\exp(\gamma \chi_{i}^{+})\exp(\gamma \chi_{k}^{+})\right)\right)\right)\right\}\nonumber\\
&\leq& d\varpi K \hat{e}_{\ominus_{\nu}
\gamma}(t_{1},t_{0})\|\psi\|_{0}.
\end{eqnarray*}
We can easily obtain some upper bound of the derivative
$|u_{i}^{\Delta}(t_{1})|$ and $|u_{i}^{\Delta}(t_{1})|$ as follow:
\begin{equation}
|u_{i}^{\Delta}(t_{1})|\leq d\varpi K \hat{e}_{\ominus_{\nu}
\gamma}(t_{1},t_{0})\|\psi\|_{0},
\end{equation}
and
\begin{equation}\label{from1}
|v_{j}^{\Delta}(t_{1})|\leq d\varpi K \hat{e}_{\ominus_{\nu}
\gamma}(t_{1},t_{0})\|\psi\|_{0}.
\end{equation}
From (\ref{from})-(\ref{from1}), we obtain
\begin{equation}
\|Z(t_{1})-Z^{*}(t_{1})\|< d\varpi K \hat{e}_{\ominus_{\nu}
\gamma}(t_{1},t_{0})\|\psi\|_{0}.
\end{equation}
which contradicts (\ref{abcd73}), therefore (\ref{abcd72}) holds.
Letting $\varpi\longrightarrow1$, then (\ref{abcd71}) holds. Which
implies that only Stepanov-like weighted pseudo-almost automorphic
on time-space scales solution of system (\ref{eq1}) is globally
exponentially stable.
\subsection{Convergence}
\begin{definition} (\cite{adnene+ahmed}) For each $t\in\mathbb{T}$, let $N$ be a neighborhood of $t$.
Then, we define the generalized derivative (or Dini derivative on
time-space scales) $D^{+}V^{\Delta}(t)$, to mean that, given
$\epsilon>0$, there exists a right neighborhood $N_{\epsilon}\subset
N$ of $t$ such
\begin{equation*}
D^{+}V^{\Delta}(t)=D^{+}V^{\Delta}(t,x(t))=\frac{V(\sigma(t)),x(\sigma(t))-V(t,x(t)}{\nu(t)}.
\end{equation*}
\end{definition}
\begin{theorem}\label{th3}
Suppose that assumptions ($H_{1}$)-($H_{4}$) hold.\\
Let $h^{\ast}\left( \cdot\right)=\left(x_{1}^{\ast
}\left(\cdot\right),\cdots
,x_{n}^{\ast}\left(\cdot\right),y_{1}^{\ast}\left(\cdot\right),...,y_{m}^{\ast}\left(\cdot\right)\right)^{T}$
be a Stepanov-like weighted pseudo-almost automorphic on time-space
scales solution of system (\ref{eq1}). If
\begin{multline*}
\alpha_{i}^{-}-\sum\limits_{j=1}^{m}\left[L_{j}\left(D_{ij}^{+}+(D_{ij}^{\tau})^{+}
+(\overline{D}_{ij})^{+}\sigma_{ij}^{+}+(\widetilde{D}_{ij})^{+}\xi_{ij}^{+}
+\sum\limits_{k=1}^{m}T_{ijk}^{+}(L_{k}r+|f_{k}(0)|)\right)\right.\\
+\left.(L_{j}r+|f_{j}(0)|)\sum\limits_{k=1}^{m}T_{ijk}^{+}L_{k}\right]>0,
\end{multline*}
and
\begin{multline*}
c_{j}^{-}-\sum\limits_{i=1}^{n}\left[L_{i}\left(E_{ij}^{+}+(E_{ij}^{\tau})^{+}
+(\overline{E}_{ij})^{+}\sigma_{ij}^{+}+(\widetilde{E}_{ij})^{+}\xi_{ij}^{+}
+\sum\limits_{k=1}^{n}T_{ijk}^{+}(L_{k}r+|f_{k}(0)|)\right)\right.\\
+\left.(L_{i}r+|f_{i}(0)|)\sum\limits_{k=1}^{n}T_{ijk}^{+}L_{k}\right]>0,
\end{multline*}
 then all solutions
$\psi=(\varphi_{1},...,\varphi_{n},\phi_{1},...,\phi_{n})$ of
(\ref{eq1}) satisfying
\begin{equation*}
x_{i}^{\ast}\left(0\right)=\varphi_{i}\left(0\right),\,\
y_{j}^{\ast}\left(0\right)=\phi_{j}\left(0\right),1\leq i\leq n, \,\
1\leq j\leq m
\end{equation*}
converge to its unique Stepanov-like weighted pseudo-almost
automorphic on time-space scales solution $h^{\ast}$.
\end{theorem}
Proof. Let
$h^{\ast}\left(\cdot\right)=\left(x_{1}^{*}\left(\cdot\right),...,x_{n}^{*}\left(\cdot\right),
y_{1}^{*}\left(\cdot\right),...,y_{m}^{*}\left(\cdot \right)\right)$
be a solution of (\ref{eq1}) and
$\psi\left(\cdot\right)=\left(\varphi_{1}\left(\cdot
\right),...,\varphi_{n}\left(\cdot\right),\phi_{1}\left(\cdot\right),...,\phi_{m}\left(\cdot\right)\right)$
be a Stepanov-like weighted pseudo almost automorphic on time-space
scales solution of (\ref{eq1}). First, one verifies without
difficulty that
\begin{eqnarray*}
&&\left(x_{i}^{\ast}\left(t\right)-\alpha_{i}(t)\int_{t-\eta_{i}(t)}^{t}x_{i}(u)\Delta
u\right)^{\Delta}
-\left(\varphi_{i}^{\ast}\left(t\right)-\alpha_{i}(t)\int_{t-\eta_{i}(t)}^{t}\varphi_{i}(u)\Delta u\right)^{\Delta}\\
&=&-\alpha_{i}(t)\left(x_{i}^{\ast}\left(t-\eta_{i}(t)\right)-\varphi_{i}\left(t-\eta_{i}(t)\right)
\right)+\sum\limits_{j=1}^{m}D_{ij}\left(t\right)\left[f_{j}(x_{j}^{\ast}\left(
t\right))-f_{j}(\varphi_{j}\left(t\right))\right]\\
&+&\sum\limits_{j=1}^{m}D_{ij}^{\tau}\left(t\right)\left[f_{j}(x_{j}^{\ast}\left(t-\tau_{ij}(t)\right))-f_{j}(\varphi_{j}\left(t-\tau_{ij}(t)\right))\right]
\\&+&\sum\limits_{j=1}^{m}\overline{D}_{ij}\left(t\right)\int\limits_{t-\sigma_{ij}(t)}^{t}\left[
f_{j}(x_{j}^{\ast}\left(t-u\right))
-f_{j}(\varphi_{j}\left( t-u\right))\right]\Delta u\\
&+&\sum\limits_{j=1}^{m}\widetilde{D}_{ij}\left(t\right)\int\limits_{t-\xi_{ij}(t)}^{t}\left(f_{j}\left((x_{j}^{\ast})^{\Delta}\left(u\right)\right)-f_{j}\left(\varphi_{j}^{\Delta}\left(u\right)\right)\right)\Delta
u\\&+&\sum\limits_{i=1}^{m}\sum\limits_{j=1}^{m}\sum\limits_{k=1}^{m}T_{ijk}(t)
\left\vert
f_{k}(x_{k}(t-\chi_{k}(t)))f_{j}(x_{j}(t-\chi_{j}(t)))-f_{k}(\varphi_{k}(t-\chi_{k}(t)))f_{j}(\varphi_{j}(t-\chi_{j}(t)))\right\vert,
\end{eqnarray*}
and
\begin{eqnarray*}
&&\left(y_{j}^{\ast}\left(t\right)-c_{j}(t)\int_{t-\varsigma_{j}(t)}^{t}y_{j}(u)\Delta
u\right)^{\Delta}
-\left(\phi_{j}^{\ast}\left(t\right)-c_{j}(t)\int_{t-\varsigma_{j}(t)}^{t}\phi_{j}(u)\Delta u\right)^{\Delta}\\
&=&-c_{j}(t)\left(y_{j}^{\ast}\left(t-\varsigma_{j}(t)\right)-\phi_{j}\left(t-\varsigma_{j}(t)\right)
\right)
+\sum\limits_{j=1}^{n}E_{ij}\left(t\right)\left[f_{j}(x_{j}^{\ast}\left(
t\right))-f_{j}(\varphi_{j}\left(t\right))\right]\\
&+&\sum\limits_{j=1}^{n}E_{ij}^{\tau}\left(t\right)\left[f_{j}(x_{j}^{\ast}\left(t-\tau_{ij}(t)\right))-f_{j}(\varphi_{j}\left(t-\tau_{ij}(t)\right))\right]
\\&+&\sum\limits_{j=1}^{n}\overline{E}_{ij}\left(t\right)\int\limits_{t-\sigma_{ij}(t)}^{t}\left[
f_{j}(x_{j}^{\ast}\left(t-u\right))
-f_{j}(\varphi_{j}\left( t-u\right))\right]\Delta u\\
&+&\sum\limits_{j=1}^{n}\widetilde{E}_{ij}\left(t\right)\int\limits_{t-\xi_{ij}(t)}^{t}\left(f_{j}\left((x_{j}^{\ast})^{\Delta}\left(u\right)\right)-f_{j}\left(\varphi_{j}^{\Delta}\left(u\right)\right)\right)\Delta u\\
&+&\sum\limits_{i=1}^{n}\sum\limits_{j=1}^{n}\sum\limits_{k=1}^{n}\overline{T}_{ijk}(t)
\left\vert
f_{k}(x_{k}(t-\chi_{k}(t)))f_{j}(x_{j}(t-\chi_{j}(t)))-f_{k}(\varphi_{k}(t-\chi_{k}(t)))f_{j}(\varphi_{j}(t-\chi_{j}(t)))\right\vert,
\end{eqnarray*}
Now, consider the following (ad-hoc) Lyapunov-Krasovskii functional
\begin{equation*}
\begin{array}{cccc}
V: & \mathbb{R} & \longrightarrow & S^{p}WPAP\left(\mathbb{T},\mathbb{R}^{n}\right) \\
& t & \longmapsto & V_{1}(t)+V_{2}(t)+V_{3}(t)+V_{4}(t)+V_{5}(t),
\end{array}
\end{equation*}
where
\begin{eqnarray*}
V_{1}(t)&=&\sum\limits_{i=1}^{n}\left\vert \left(x_{i}^{\ast }\left(
t\right)-\alpha_{i}(t)\int_{t-\eta_{i}(t)}^{t}x_{i}(u)\Delta
u\right)
-\left(\varphi_{i}^{\ast }\left( t\right)-\alpha_{i}(t)\int_{t-\eta_{i}(t)}^{t}\varphi_{i}(u)\Delta u\right)\right\vert,\\
V_{2}\left(t\right)&=&\sum\limits_{j=1}^{n}\sum\limits_{i=1}^{n}\int\limits_{t-\tau_{j}(t)
}^{t}L_{j}\left(D_{ij}^{+}+(D_{ij}^{\tau})^{+}\right)\left\vert
x_{i}^{\ast}\left( s\right)
-\varphi_{i}\left( s\right)\right\vert\Delta s,\\
V_{3}\left(t\right)
&=&\sum\limits_{j=1}^{n}\sum\limits_{i=1}^{n}\int\limits_{t-\sigma_{ij}(t)}^{t
}\int\limits_{t-s}^{t}L_{j}\overline{D}_{ij}^{+}\left\vert
x_{i}^{\ast
}\left(u\right)-\varphi_{i}\left( u\right)\right\vert \Delta u \Delta s,\\
V_{4}\left( t\right)
&=&\sum\limits_{j=1}^{n}\sum\limits_{i=1}^{n}\int\limits_{t-\xi_{ij}(t)}^{t}\int\limits_{s}^{0}L_{j}\widetilde{D}_{ij}^{+}\left\vert
(x_{i}^{\ast})^{\Delta}\left(
u\right)-(\varphi_{i})^{\Delta}\left(u\right)\right\vert\Delta
u\Delta s,
\end{eqnarray*}
and
\begin{eqnarray*}
V_{5}(t)&=&\sum\limits_{j=1}^{n}\left\vert \left(y_{j}^{\ast }\left(
t\right)-c_{j}(t)\int_{t-\varsigma_{j}(t)}^{t}y_{j}(u)\Delta
u\right)
-\left(\phi_{j}^{\ast}\left(t\right)-c_{j}(t)\int_{t-\varsigma_{j}(t)}^{t}\phi_{j}(u)\Delta
u\right)\right\vert.
\end{eqnarray*}
Let us calculate the upper right Dini derivative on time-space
scales $D^{+}V^{\Delta}\left(t\right)$ of $V$ along the trajectory
of the solution of the equation above. Then one has
\begin{eqnarray*}
D^{+}V_{1}^{\Delta}(t) &\leq
&-\sum\limits_{i=1}^{m}\alpha_{i}^{-}\left\vert
x_{i}^{\ast}\left(t\right)-\varphi_{i}\left(t\right)\right\vert
+\sum\limits_{i=1}^{m}\sum\limits_{j=1}^{m}D_{ij}^{+} L_{j}\left\vert x_{j}^{\ast}\left(t\right)-\varphi _{j}\left( t \right) \right\vert\\
&+&\sum\limits_{i=1}^{m}\sum\limits_{j=1}^{m}(D_{ij}^{\tau})^{+}L_{j}\left\vert
x_{j}^{\ast}\left(t-\tau_{j}(t)\right)
-\varphi_{j}\left(t-\tau_{j}(t)\right)\right\vert \\
&+&\sum\limits_{i=1}^{m}\sum\limits_{j=1}^{m}\overline{D}_{ij}^{+}\int\limits_{t-\sigma_{ij}(t)}^{t}L_{j}\left\vert(x_{j}^{\ast}\left(
t-u\right)
-\varphi_{j}\left( t-u\right) \right\vert \Delta u\\
&+&\sum\limits_{i=1}^{m}\sum\limits_{j=1}^{m}
\widetilde{D}_{ij}^{+}\int\limits_{t-\xi_{ij}(t)}^{t}
L_{j}\left\vert(x_{j}^{\ast})^{\Delta}
\left( u\right)-\varphi_{j}^{\Delta }\left( u\right)\right\vert \Delta u\\
&+&\sum\limits_{i=1}^{m}\sum\limits_{j=1}^{m}\sum\limits_{k=1}^{m}T_{ijk}^{+} (L_{k}r+|f_{k}(0)|)L_{j}\left\vert x_{j}(t-\chi_{j}(t))-\varphi_{j}(t-\chi_{j}(t))\right\vert\\
&+&\sum\limits_{i=1}^{m}\sum\limits_{j=1}^{m}\sum\limits_{k=1}^{m}T_{ijk}^{+}(L_{j}r+|f_{j}(0)|)L_{k}\left\vert
x_{k}(t-\chi_{k}(t))-\varphi_{k}(t-\chi_{k}(t))\right\vert
\end{eqnarray*}
Obviously,
\begin{eqnarray*}
D^{+}V_{2}^{\Delta}(t) &\leq
&\sum\limits_{j=1}^{m}\sum\limits_{i=1}^{m}L_{j}(D_{ij}^{+}+(D_{ij}^{\tau})^{+})\left[\left\vert
x_{i}^{\ast }\left(t\right)-\varphi_{i}\left(t\right)
\right\vert-\left\vert x_{i}^{\ast }\left(t-\tau_{ij}(t) \right)-\varphi_{i}\left(t-\tau_{ij}(t)\right) \right\vert \right]\\
&\leq&\sum\limits_{i=1}^{m}\sum\limits_{j=1}^{m}L_{j}(D_{ij}^{+}
+(D_{ij}^{\tau})^{+})\left\vert x_{i}^{\ast}\left( t\right)-\varphi_{i}\left( t\right) \right\vert\\
&-&\sum\limits_{i=1}^{m}\sum\limits_{j=1}^{m}L_{j}(D_{ij}^{+}+(D_{ij}^{\tau})^{+})\left\vert
x_{i}^{\ast}\left(t-\tau_{j}(t)\right)-\varphi_{i}\left(t-\tau_{j}(t)\right)\right\vert
\end{eqnarray*}
and
\begin{eqnarray*}
D^{+}V_{3}^{\Delta}(t)&\leq&\sum\limits_{j=1}^{m}\sum\limits_{i=1}^{m}L_{j}\overline{D}_{ij}^{+}\int_{t-\sigma_{ij}(t)}^{t}
\left[\left\vert x_{i}^{\ast }\left(t\right)-\varphi_{i}\left(
t\right) \right\vert
-\left\vert x_{i}^{\ast}\left(t-s\right)-\varphi_{i}\left(t-s\right)\right\vert\right]\Delta s \\
&\leq&\sum\limits_{j=1}^{m}\sum\limits_{i=1}^{m}L_{j}\overline{D}_{ij}^{+}\int_{t-\sigma_{ij}(t)}^{t}\left\vert
x_{i}^{\ast }\left( t\right)-\varphi_{i}\left( t\right) \right\vert \Delta s\\
&-&\sum\limits_{j=1}^{m}\sum\limits_{i=1}^{m}L_{j}
\overline{D}_{ij}^{+}\int_{t-\sigma_{ij}(t)}^{t}\left\vert
x_{i}^{\ast }\left(t-s\right)-\varphi_{i}\left(t-s\right)\right\vert
\Delta s.
\end{eqnarray*}
Reasoning in a similar way, we can obtain the following estimation
\begin{eqnarray*}
D^{+}V_{4}^{\Delta}\left(t\right)&=&\sum\limits_{j=1}^{m}\sum\limits_{i=1}^{m}\int_{t-\xi_{ij}(t)}^{t}\int\limits_{s}^{0}L_{j}\widetilde{D}_{ij}^{+}\left\vert
(x_{i}^{\ast})^{\Delta}\left(u\right)-\varphi_{i}^{\Delta}\left(u\right)
\right\vert \Delta u\,\ \Delta s \\
&\leq&-\sum\limits_{j=1}^{m}\sum\limits_{i=1}^{m}L_{j}\widetilde{D}_{ij}^{+}\int_{t-\xi_{ij}(t)}^{t}\left\vert
(x_{i}^{\ast})^{\Delta}\left(s\right)-\varphi_{i}^{\Delta}\left(s\right)\right\vert
\Delta s,
\end{eqnarray*}
and,
\begin{eqnarray*}
D^{+}V_{5}^{\Delta}\left(t\right)&\leq
&-\sum\limits_{j=1}^{n}c_{j}^{-}\left\vert
y_{j}^{\ast}\left(t\right)-\phi_{j}\left(t\right)\right\vert
+\sum\limits_{i=1}^{n}\sum\limits_{j=1}^{n}E_{ij}^{+} L_{j}\left\vert y_{j}^{\ast}\left(t\right)-\phi _{j}\left( t \right) \right\vert\\
&+&\sum\limits_{i=1}^{n}\sum\limits_{j=1}^{n}(E_{ij}^{\tau})^{+}L_{j}\left\vert
y_{j}^{\ast}\left(t-\tau_{j}(t)\right)
-\phi_{j}\left(t-\tau_{j}(t)\right)\right\vert \\
&+&\sum\limits_{i=1}^{n}\sum\limits_{j=1}^{n}\overline{E}_{ij}^{+}\int\limits_{t-\sigma_{ij}(t)}^{t}L_{j}\left\vert(y_{j}^{\ast}\left(
t-u\right)
-\phi_{j}\left( t-u\right) \right\vert \Delta u\\
&+&\sum\limits_{i=1}^{n}\sum\limits_{j=1}^{n}
\widetilde{E}_{ij}^{+}\int\limits_{t-\xi_{ij}(t)}^{t}
L_{j}\left\vert(y_{j}^{\ast})^{\Delta} \left(
u\right)-\phi_{j}^{\Delta }\left( u\right)\right\vert \Delta
u+\left\vert y_{j}(t)-\phi_{j}(t)\right\vert
\end{eqnarray*}
\begin{eqnarray*}
&+&\sum\limits_{i=1}^{n}\sum\limits_{j=1}^{n}\sum\limits_{k=1}^{n}\overline{T}_{ijk}^{+} (L_{k}r+|f_{k}(0)|)L_{j}\left\vert y_{j}(t-\chi_{j}(t))-\phi_{j}(t-\chi_{j}(t))\right\vert\\
&+&\sum\limits_{i=1}^{n}\sum\limits_{j=1}^{n}\sum\limits_{k=1}^{n}\overline{T}_{ijk}^{+}(L_{j}r+|f_{j}(0)|)L_{k}\left\vert
y_{k}(t-\chi_{k}(t))-\varphi_{k}(t-\chi_{k}(t))\right\vert
\end{eqnarray*}
By using the inequality of the Dini derivative on time-space scales
\begin{equation*}
D^{+}\left(F_{1}^{\Delta}+F_{2}^{\Delta}\right)\leq
D^{+}\left(F_{1}^{\Delta}\right) +D^{+}\left(F_{2}^{\Delta}\right),
\end{equation*}
we get
\begin{eqnarray*}
D^{+}\left( V^{\Delta}(t)\right) &\leq
&D^{+}V_{1}^{\Delta}(t)+D^{+}V_{2}^{\Delta}(t)+D^{+}V_{3}^{\Delta}(t)+D^{+}V_{4}^{\Delta}(t) \\
&\leq &-\sum\limits_{i=1}^{m}
\sum\limits_{j=1}^{n}\min\left\{\alpha_{i}^{+}-D_{ij}^{+}L_{j}-(D_{ij}^{\tau})^{+}L_{j}-\overline{D}_{ij}^{+}\sigma_{ij}^{+}L_{j}\right.\\
&-&\widetilde{D}_{ij}^{+}\xi_{ij}^{+}L_{j}+L_{j}\sum\limits_{k=1}^{n}T_{ijk}^{+}(L_{k}r+|f_{k}(0)|)
+(L_{j}r+|f_{j}(0)|)\sum\limits_{k=1}^{n}T_{ijk}^{+}L_{k},\\
&&
c_{j}^{+}-E_{ij}^{+}L_{j}-(E_{ij}^{\tau})^{+}L_{j}-\overline{E}_{ij}^{+}\sigma_{ij}^{+}L_{j}
-\widetilde{E}_{ij}^{+}\xi_{ij}^{+}L_{j}+L_{j}\sum\limits_{k=1}^{m}\overline{T}_{ijk}^{+}(L_{k}r+|f_{k}(0)|)\\
&+&\left.(L_{j}r+|f_{j}(0)|)\sum\limits_{k=1}^{m}\overline{T}_{ijk}^{+}L_{k}\right\} \| h^{\ast}\left(t\right) -\psi\left( t\right) \| \\
&=&-\sum\limits_{i=1}^{n}\sum\limits_{j=1}^{m}\min\{\beta
_{i},\beta_{j}\}\| h^{\ast }\left(t\right)-\psi\left(t\right)\| <0.
\end{eqnarray*}
By integrating the above inequality from $t_{0}$ to $t$, we get
\begin{equation*}
V(t)+\sum\limits_{i=1}^{n}\sum\limits_{j=1}^{m}\min\{\beta
_{i},\beta_{j}\}\int\limits_{t_{0}}^{t}\|
y^{\ast}\left(t\right)-\psi\left(t\right)\|\Delta
s<V(t_{0})<+\infty.
\end{equation*}
Now, we remark that $V(t)>0.$ It follows that
\begin{equation*}
\lim\limits_{t\rightarrow+\infty}\sup\int\limits_{t_{0}}^{t}\min\{\beta
_{i},\beta_{j}\}\| h^{\ast}\left(s\right) -\psi\left( s\right)\|
\Delta s<V(t_{0})<+\infty.
\end{equation*}
Note that $h^{\ast }\left(\cdot\right)$ is bounded on
$\mathbb{T}^{+}$. Therefore
\begin{equation*}
\lim\limits_{t\rightarrow+\infty}\|
h^{\ast}\left(t\right)-\psi\left(t\right)\|=0.
\end{equation*}
The proof of this theorem is now completed.
\begin{remark}
Theorem \ref{th0}, Theorem \ref{th1}, Theorem \ref{th2} and Theorem
\ref{th3} are new even for the both cases of differential equations
($\mathbb{T}=\mathbb{R}$) and difference equations
($\mathbb{T}=\mathbb{Z}$).
\end{remark}
\section{Numerical example}
In system (\ref{eq1}), let $n=3,\,\ m=2$, and take coefficients as
follows:
\begin{multline*}
\left(
\begin{array}{ccc}
D_{11}(t) &  D_{12}(t) \\
D_{21}(t) & D_{22}(t) \\
D_{31}(t) &  D_{32}(t) \\
\end{array}
\right)=\left(
\begin{array}{ccc}
D^{\tau}_{11}(t) &  D^{\tau}_{12}(t) \\
D^{\tau}_{21}(t) & D^{\tau}_{22}(t) \\
D^{\tau}_{31}(t) &  D^{\tau}_{32}(t) \\
\end{array}
\right)=\left(
\begin{array}{ccc}
\tilde{D}_{11}(t) &  \tilde{D}_{12}(t) \\
\tilde{D}_{21}(t) & \tilde{D}_{22}(t) \\
\tilde{D}_{31}(t) &  \tilde{D}_{32}(t) \\
\end{array}
\right)\\=\left(
\begin{array}{ccc}
\bar{D}_{11}(t) &  \bar{D}_{12}(t) \\
\bar{D}_{21}(t) & \bar{D}_{22}(t) \\
\bar{D}_{31}(t) &  \bar{D}_{32}(t) \\
\end{array}
\right)=\left(
\begin{array}{ccc}
\frac{\sin(t)}{20} &  \frac{\sin(t)}{20} \\
\frac{\cos(t)}{20} &  \frac{\sin(t)}{20} \\
\frac{\cos(t)}{20} &  \frac{\sin(t)}{20} \\
\end{array}
\right);
\\
\left(
\begin{array}{ccc}
E_{11}(t) &  E_{12}(t) \\
E_{21}(t) & E_{22}(t) \\
E_{31}(t) &  E_{32}(t) \\
\end{array}
\right)=\left(
\begin{array}{ccc}
E^{\tau}_{11}(t) &  E^{\tau}_{12}(t) \\
E^{\tau}_{21}(t) & E^{\tau}_{22}(t) \\
E^{\tau}_{31}(t) &  E^{\tau}_{32}(t) \\
\end{array}
\right)=\left(
\begin{array}{ccc}
\tilde{E}_{11}(t) &  \tilde{E}_{12}(t) \\
\tilde{E}_{21}(t) & \tilde{E}_{22}(t) \\
\tilde{E}_{31}(t) &  \tilde{E}_{32}(t) \\
\end{array}
\right)
\\=\left(
\begin{array}{ccc}
\bar{E}_{11}(t) &  \bar{E}_{12}(t) \\
\bar{E}_{21}(t) & \bar{E}_{22}(t) \\
\bar{E}_{31}(t) &  \bar{E}_{32}(t) \\
\end{array}
\right)=\left(
\begin{array}{ccc}
\frac{\sin(t)}{15} &  \frac{\sin(t)}{15} \\
\frac{\cos(t)}{15} &  \frac{\sin(t)}{15} \\
\frac{\cos(t)}{15} &  \frac{\sin(t)}{15} \\
\end{array}
\right);\\
\alpha_{1}(t)=\alpha_{2}(t)=0.73+0.02\sin\frac{1}{2+\cos(t)+\cos(\sqrt{2}t)};\\
f_{1}(x)=f_{2}(x)=0.1\arctan x;\,\ \nu(t)=\hat{e}(t,\sigma(t)) \,\ \text{for all} \,\ t\in[0,\infty)_{\mathbb{T}}; \\ \nu(t)=1 \,\ \text{for all}\,\ t\in(-\infty,0)_{\mathbb{T}};\\
c_{1}(t)=0.54-0.02\sin\frac{1}{2+\cos(t)+\cos(\sqrt{2}t)};\\
c_{2}(t)=0.54+0.02\sin\frac{1}{2+\cos(t)+\cos(\sqrt{2}t)};\\
I_{1}(t)=0.01\sin\frac{1}{2+\sin(t)+\sin(\sqrt{2}t)}+0.01e^{-t^{4}\cos^{2}(t)};\\
I_{2}(t)=0.02\sin\frac{1}{2+\sin(t)+\sin(\sqrt{2}t)}+0.02e^{-t^{4}\cos^{2}(t)};\phantom{+++++}\\
J_{1}(t)=0.02\sin\frac{1}{2+\cos(t)+\cos(\sqrt{2}t)}+0.02e^{-t^{2}\sin^{4}(t)};\\
J_{2}(t)=0.1\sin\frac{1}{2+\sin(t)+\sin(\sqrt{2}t)};\\
\eta_{1}(t)=0.04\sin\frac{1}{2+\sin(t)+\sin(\sqrt{2}t)}+0.03e^{-t^{4}\cos^{2}(t)};\\
\eta_{2}(t)=0.01\sin\frac{1}{2+\cos(t)+\cos(\sqrt{2}t)}+0.01e^{-t^{4}\cos^{2}(t)};\\
\varsigma_{1}(t)=0.01\sin\frac{1}{2+\cos(t)+\cos(\sqrt{2}t)}+0.01e^{-t^{2}\sin^{4}(t)};\\
\varsigma_{2}(t)=0.01\sin\frac{1}{2+\cos(t)+\cos(\sqrt{2}t)}+0.01e^{-t^{2}\sin^{4}(t)};\\
\tau_{11}(t)=0.02\sin(\sqrt{2}t)+e^{-t^{2}\sin^{4}(t)};\\
\tau_{12}(t)=0.02\sin\frac{1}{2+\sin(t)+\sin(\sqrt{2}t)}+0.01e^{-t^{2}\sin^{4}(t)};\\
\tau_{21}(t)=0.01\sin\frac{1}{2+\cos(t)+\cos(\sqrt{2}t)}+0.01e^{-t^{4}\cos^{2}(t)};\\
\tau_{22}(t)=0.02\sin\frac{1}{2+\sin(t)+\sin(\sqrt{2}t)}+0.01e^{-t^{4}\cos^{2}(t)};\phantom{++++}
\end{multline*}
\begin{eqnarray*}
\sigma_{11}(t)&=&0.02\sin\frac{1}{2+\cos(t)+\cos(\sqrt{2}t)}+0.02e^{-t^{2}\sin^{4}(t)};\\
\sigma_{12}(t)&=&0.01\sin\frac{1}{2+\cos(t)+\cos(\sqrt{2}t)}+0.01e^{-t^{2}\sin^{4}(t)};\\
\sigma_{21}(t)&=&0.02\sin\frac{1}{2+\sin(t)+\sin(\sqrt{2}t)}+0.01e^{-t^{4}\cos^{2}(t)};\\
\sigma_{12}(t)&=&0.02\sin\frac{1}{2+\cos(t)+\cos(\sqrt{2}t)}+e^{-t^{2}\sin^{4}(t)};\\
\xi_{11}(t)&=&0.02\sin\frac{1}{2+\cos(t)+\cos(\sqrt{2}t)}+0.01e^{-t^{4}\cos^{2}(t)};\\
\xi_{12}(t)&=&0.03\sin\frac{1}{2+\sin(t)+\sin(\sqrt{2}t)}+0.02e^{-t^{2}\sin^{4}(t)};\\
\xi_{21}(t)&=&0.01\sin\frac{1}{2+\cos(t)+\cos(\sqrt{2}t)}+0.01e^{-t^{4}\cos^{2}(t)};\\
\xi_{22}(t)&=&0.02\sin\frac{1}{2+\sin(t)+\sin(\sqrt{2}t)}+0.02e^{-t^{2}\sin^{4}(t)}.
\end{eqnarray*}
By a simple calculation, we have
\begin{multline*}
\left(
\begin{array}{ccc}
D_{11}^{+} &  D_{12}^{+} \\
D_{21}^{+} & D_{22}^{+} \\
D_{31}^{+} &  D_{32}^{+} \\
\end{array}
\right)=\left(
\begin{array}{ccc}
(D^{\tau}_{11})^{+} &  (D^{\tau}_{12})^{+} \\
(D^{\tau}_{21})^{+} & (D^{\tau}_{22})^{+} \\
(D^{\tau}_{31})^{+} &  (D^{\tau}_{32})^{+} \\
\end{array}
\right)=\left(
\begin{array}{ccc}
\tilde{D}_{11}^{+} &  \tilde{D}_{12}^{+}\\
\tilde{D}_{21}^{+} & \tilde{D}_{22}^{+} \\
\tilde{D}_{31}^{+} &  \tilde{D}_{32}^{+}\\
\end{array}
\right)=\left(
\begin{array}{ccc}
\bar{D}_{11}^{+} &  \bar{D}_{12}^{+} \\
\bar{D}_{21}^{+} & \bar{D}_{22}^{+} \\
\bar{D}_{31}^{+} &  \bar{D}_{32}^{+} \\
\end{array}
\right)=\left(
\begin{array}{ccc}
\frac{1}{20} &  \frac{1}{20} \\
\frac{1}{20} & \frac{1}{20} \\
\frac{1}{20} &  \frac{1}{20} \\
\end{array}
\right);
\\
\left(
\begin{array}{ccc}
E_{11}^{+} &  E_{12}^{+} \\
E_{21}^{+} & E_{22}^{+} \\
E_{31}^{+} &  E_{32}^{+} \\
\end{array}
\right)=\left(
\begin{array}{ccc}
(E^{\tau}_{11})^{+} &  (E^{\tau}_{12})^{+} \\
(E^{\tau}_{21})^{+} & (E^{\tau}_{22})^{+} \\
(E^{\tau}_{31})^{+} &  (E^{\tau}_{32})^{+} \\
\end{array}
\right)=\left(
\begin{array}{ccc}
\tilde{E}_{11}^{+} &  \tilde{E}_{12}^{+}\\
\tilde{E}_{21}^{+} & \tilde{E}_{22}^{+} \\
\tilde{E}_{31}^{+} &  \tilde{E}_{32}^{+}\\
\end{array}
\right)=\left(
\begin{array}{ccc}
\bar{E}_{11}^{+} &  \bar{E}_{12}^{+} \\
\bar{E}_{21}^{+} & \bar{E}_{22}^{+} \\
\bar{E}_{31}^{+} &  \bar{E}_{32}^{+} \\
\end{array}
\right)=\left(
\begin{array}{ccc}
\frac{1}{15} &  \frac{1}{15} \\
\frac{1}{15} & \frac{1}{15} \\
\frac{1}{15} &  \frac{1}{15} \\
\end{array}
\right);\\
I_{1}^{+}=0.02; \,\ I_{2}^{+}=0.04;\,\ I_{3}^{+}=0.05;\,\ J_{1}^{+}=0.01;\,\ J_{2}^{+}=-0.1;\\
\eta_{1}^{+}=0.07; \,\ \eta_{2}^{+}=0.02;\,\ \eta_{3}^{+}= \,\ \varsigma_{1}^{+}=0.02; \varsigma_{2}^{+}=0.02;\\
(\sigma_{ij}^{+})_{1\leq i\leq 3;1\leq j\leq 2}=\left(
\begin{array}{ccc}
1.04 & 1.02 \\
1.03 & 1.02 \\
1.03 & 1.02 \\
\end{array}
\right); (\xi_{ij}^{+})_{1\leq i\leq 3;1\leq j\leq 2}=\left(
\begin{array}{ccc}
1.03 & 1.05 \\
1.02 & 1.04 \\
1.02 & 1.04 \\
\end{array}
\right);\phantom{+++++++++++++}\\
(T_{1jk})_{2\times2}=\left[\begin{array}{ccc}
0.01\sin \sqrt{2}t & 0.01\cos \sqrt{2}t \\
0.01\cos \sqrt{2}t & 0.01\frac{\cos \sqrt{2}t}{2} \\
\end{array}
\right];\,\ (T_{2jk})_{2\times2}=\left[
\begin{array}{ccc}
0.01\sin \sqrt{2}t & 0.01\frac{\cos \sqrt{2}t}{2}  \\
0.01\sin \sqrt{2}t & 0.01\sin \sqrt{2}t  \\
\end{array}
\right];\\ (T_{3jk})_{2\times2}=\left[
\begin{array}{ccc}
0.01\sin \sqrt{2}t & 0.01\frac{\cos \sqrt{2}t}{2}  \\
0.01\sin \sqrt{2}t & 0.01\sin \sqrt{2}t  \\
\end{array}
\right]; (\overline{T}_{1jk})_{2\times2}=\left[\begin{array}{ccc}
0.01\sin \sqrt{2}t & 0.01\cos \sqrt{2}t1 \\
0.01\cos \sqrt{2}t & 0.01\frac{\cos \sqrt{2}t}{2} \\
\end{array}
\right];\\ (\overline{T}_{2jk})_{2\times2}=\left[
\begin{array}{ccc}
0.01\sin \sqrt{2}t & 0.01\frac{\cos \sqrt{2}t}{2}  \\
0.01\sin \sqrt{2}t & 0.01\sin \sqrt{2}t  \\
\end{array}
\right];\,\ (\overline{T}_{3jk})_{2\times2}=\left[
\begin{array}{ccc}
0.01\sin \sqrt{2}t & 0.01\frac{\cos \sqrt{2}t}{2}  \\
0.01\sin \sqrt{2}t & 0.01\sin \sqrt{2}t  \\
\end{array}
\right];
\end{multline*}
We can take $L_{1}=L_{2}=0.1$, $r=0.43$ and we have
\begin{eqnarray*}
M_{1}&=&\alpha_{1}^{+}\eta_{1}^{+}r+\sum\limits_{j=1}^{2}\left(D_{1j}^{+}+(D_{1j}^{\tau})^{+}+\overline{D}_{1j}^{+}\sigma_{1j}^{+}\right.
+\left.\widetilde{D}_{1j}^{+}\xi_{1j}^{+}\right)(L_{j}r+|f_{j}(0)|)\\
&+&\sum\limits_{j=1}^{2}\sum\limits_{k=1}^{2}T_{1jk}^{+}(L_{k}r+|f_{k}(0)|)(L_{j}r+|f_{j}(0)|)+I_{1}^{+},\\
&=& 0.75\times0.07\times0.43+\left(\frac{1}{15}+\frac{1}{15}+\frac{1}{15}\times 1.04+\frac{1}{15}\times 1.03\right)\left(0.1\times0.43+0.1\right)\\
&+& \left(\frac{1}{15}+\frac{1}{15}+\frac{1}{15}\times 1.02+\frac{1}{15}\times 1.05\right)\left(0.1\times0.43+0.1\right)\\
&+&0.04\left(0.1\times0.43+0.1\right)^{2}+0.02\\
&=& 0.119\\
\overline{M}_{1}&=&\alpha_{1}^{+}\eta_{1}^{+}+\sum\limits_{j=1}^{2}\left(D_{1j}^{+}+(D_{1j}^{\tau})^{+}+\overline{D}_{1j}^{+}\sigma_{ij}^{+}
+\widetilde{D}_{1j}^{+}\xi_{1j}^{+}\right)L_{j}\\
&+& \sum\limits_{j=1}^{2}\sum\limits_{k=1}^{2}(T_{1jk}^{+}+T_{1kj}^{+})(L_{k}r+|f_{k}(0)|)\\
&=& 0.75\times0.07+\left(\frac{1}{15}+\frac{1}{15}+\frac{1}{15}\times 1.03+\frac{1}{15}\times 1.02\right)\left(0.1\times0.43+0.1\right)\\
&+&\left(\frac{1}{15}+\frac{1}{15}+\frac{1}{15}\times 1.02+\frac{1}{15}\times 1.04\right)\left(0.1\times0.43+0.1\right)+0.08\left(0.1\times0.43+0.1\right)\\
&=& 0.139\\
N_{1}&=&c_{1}^{+}\varsigma_{1}^{+}r+\sum\limits_{i=1}^{3}\left(E_{i1}^{+}+(E_{i1}^{\tau})^{+}+\overline{E}_{i1}^{+}\sigma_{i1}^{+}
+\widetilde{E}_{i1}^{+}\xi_{i1}^{+}\right)(L_{i}r+|f_{i}(0)|)\\
&+&\sum\limits_{i=1}^{3}\sum\limits_{k=1}^{3}\overline{T}_{i1k}^{+}\left(L_{k}r+\left\vert f_{k}(0)\right\vert\right)\left(L_{i}r+\left\vert f_{i}(0)\right\vert\right)+J_{1}^{+}\\
&=& 0.56\times1.02\times0.43+\left(\frac{1}{20}+\frac{1}{20}+\frac{1}{20}\times 1.04+\frac{1}{20}\times 1.03\right)\left(0.1\times0.43+0.1\right)\\
&+& \left(\frac{1}{20}+\frac{1}{20}+\frac{1}{20}\times 1.03+\frac{1}{20}\times 1.02\right)\left(0.1\times0.43+0.1\right)+0.04\left(0.1\times0.43+0.1\right)^{2}+0.04\\
&=&0.343\\
\overline{N}_{1}&=&
c_{1}^{+}\varsigma_{1}^{+}+\sum\limits_{i=1}^{3}\left(E_{i1}^{+}+(E_{i1}^{\tau})^{+}
+\overline{E}_{i1}^{+}\sigma_{i1}^{+}+\widetilde{E}_{i1}^{+}\xi_{i1}^{+}\right)L_{i}\\
&+&\sum\limits_{i=1}^{3}\sum\limits_{k=1}^{3}(\overline{T}_{i1k}^{+}+\overline{T}_{ik1}^{+})(L_{k}r+|f_{k}(0)|)\\
&=&0.56\times 1.02+\left(\frac{6}{20}+1.04\times\frac{1}{20}+3.09\times\frac{1}{20}+2.04\times\frac{1}{20}\right)\times 0.1+0.018\\
&=&0.65 \\
M_{2}&=&\alpha_{2}^{+}\eta_{2}^{+}r+\sum\limits_{j=1}^{2}\left(D_{2j}^{+}+(D_{2j}^{\tau})^{+}+\overline{D}_{2j}^{+}\sigma_{2j}^{+}
+\widetilde{D}_{2j}^{+}\xi_{2j}^{+}\right)(L_{j}r+|f_{j}(0)|)\\
&+& \sum\limits_{j=1}^{2}\sum\limits_{k=1}^{2}(T_{2jk}^{+}+T_{2kj}^{+})(L_{k}r+|f_{k}(0)|)+I_{2}^{+}\\
&=&0.75\times0.02\times0.43+\left(\frac{1}{15}+\frac{1}{15}+\frac{1}{15}\times 1.03+\frac{1}{15}\times 1.02\right)\left(0.1\times0.43+0.1\right)\\
&+& \left(\frac{1}{15}+\frac{1}{15}+\frac{1}{15}\times 1.02+\frac{1}{15}\times 1.04\right)\left(0.1\times0.43+0.1\right)\\
&+&0.04+0.08\left(0.1\times0.43+0.1\right)\\
&=& 0.52
\end{eqnarray*}
\begin{eqnarray*}
\overline{M}_{2}&=&\alpha_{2}^{+}\eta_{2}^{+}+\sum\limits_{j=1}^{2}\left(D_{2j}^{+}+(D_{2j}^{\tau})^{+}+\overline{D}_{2j}^{+}\sigma_{2j}^{+}
+\widetilde{D}_{2j}^{+}\xi_{2j}^{+}\right)L_{j}\\
&+& \sum\limits_{j=1}^{2}\sum\limits_{k=1}^{2}T_{2jk}^{+}(L_{k}r+|f_{k}(0)|)(L_{j}r+|f_{j}(0)|)\\
&=& 0.75\times0.02+\left(\frac{1}{15}+\frac{1}{15}+\frac{1}{15}\times 1.03+\frac{1}{15}\times 1.02\right)\times0.1\\
&+& \left(\frac{1}{15}+\frac{1}{15}+\frac{1}{15}\times 1.02+\frac{1}{15}\times 1.04\right)\times0.1+0.04\left(0.1\times0.43+0.1\right)^{2}\\
&=& 0.069\\
N_{2}&=&c_{2}^{+}\varsigma_{2}^{+}r+\sum\limits_{i=1}^{3}\left(E_{i2}^{+}+(E_{i2}^{\tau})^{+}+\overline{E}_{i2}^{+}\sigma_{i2}^{+}
+\widetilde{E}_{i2}^{+}\xi_{i2}^{+}\right)(L_{i}r+|f_{i}(0)|)\\
&+&\sum\limits_{i=1}^{3}\sum\limits_{k=1}^{3}\overline{T}_{i2k}^{+}\left(L_{k}r+\left\vert f_{k}(0)\right\vert\right)\left(L_{i}r+\left\vert f_{i}(0)\right\vert\right)+J_{2}^{+}\\
&=& 0.213\\
\overline{N}_{2}&=&
c_{2}^{+}\varsigma_{2}^{+}+\sum\limits_{i=1}^{3}\left(E_{i2}^{+}+(E_{i2}^{\tau})^{+}
+\overline{E}_{i2}^{+}\sigma_{i2}^{+}+\widetilde{E}_{i2}^{+}\xi_{i2}^{+}\right)L_{i}\\
&+&\sum\limits_{i=1}^{3}\sum\limits_{k=1}^{3}(\overline{T}_{i2k}^{+}+\overline{T}_{ik2}^{+})(L_{k}r+|f_{k}(0)|)\\
&=&0.343\\
M_{3}&=&\alpha_{3}^{+}\eta_{3}^{+}r+\sum\limits_{j=1}^{2}\left(D_{3j}^{+}+(D_{3j}^{\tau})^{+}+\overline{D}_{3j}^{+}\sigma_{3j}^{+}\right.
+\left.\widetilde{D}_{3j}^{+}\xi_{3j}^{+}\right)(L_{j}r+|f_{j}(0)|)\\
&+&\sum\limits_{j=1}^{2}\sum\limits_{k=1}^{2}T_{3jk}^{+}(L_{k}r+|f_{k}(0)|)(L_{j}r+|f_{j}(0)|)+I_{3}^{+},\\
&=& 0.75\times0.07\times0.43+\left(\frac{1}{15}+\frac{1}{15}+\frac{1}{15}\times 0.04+\frac{1}{15}\times 0.03\right)\left(0.1\times0.43+0.1\right)\\
&+& \left(\frac{1}{15}+\frac{1}{15}+\frac{1}{15}\times 0.02+\frac{1}{15}\times 0.05\right)\left(0.1\times0.43+0.1\right)\\
&+&0.04\left(0.1\times0.43+0.1\right)^{2}+0.02\\
&=& 0.12\\
\overline{M}_{3}&=&\alpha_{3}^{+}\eta_{3}^{+}+\sum\limits_{j=1}^{2}\left(D_{3j}^{+}+(D_{3j}^{\tau})^{+}+\overline{D}_{3j}^{+}\sigma_{ij}^{+}
+\widetilde{D}_{3j}^{+}\xi_{3j}^{+}\right)L_{j}\\
&+& \sum\limits_{j=1}^{2}\sum\limits_{k=1}^{2}(T_{3jk}^{+}+T_{3kj}^{+})(L_{k}r+|f_{k}(0)|)\\
&=& 0.75\times0.07+\left(\frac{1}{15}+\frac{1}{15}+\frac{1}{15}\times 0.03+\frac{1}{15}\times 0.02\right)\left(0.1\times0.43+0.1\right)\\
&+&\left(\frac{1}{15}+\frac{1}{15}+\frac{1}{15}\times 0.02+\frac{1}{15}\times 0.04\right)\left(0.1\times0.43+0.1\right)+0.08\left(0.1\times0.43+0.1\right)\\
&=& 0.0136
\end{eqnarray*}
%%%%%%%%%%%%%%%%%%%%%%%%%%%%%%%%%%%%%%%%%%%%%%%%%%%%%%%%%%%%%%%%%%%%%%%%%%%%%%%%%%%%
The conditions $(H_{1})$, $(H_{2})$ and $(H_{4})$ are satisfied and
it is easy to verify
\begin{multline*}
\max\left\{\frac{M_{1}}{\alpha_{1}^{-}},\frac{M_{2}}{\alpha_{2}^{-}},\frac{M_{3}}{\alpha_{3}^{-}},
\left(1+\frac{\alpha_{1}^{+}}{\alpha_{1}^{-}}\right)M_{1},
\left(1+\frac{\alpha_{2}^{+}}{\alpha_{2}^{-}}\right)M_{2},\left(1+\frac{\alpha_{3}^{+}}{\alpha_{3}^{-}}\right)M_{3},\right.\\
\left.\frac{N_{1}}{c_{1}^{-}},\left(1+\frac{c_{1}^{+}}{c_{1}^{-}}\right)N_{1},
\frac{N_{2}}{c_{2}^{-}},\left(1+\frac{c_{2}^{+}}{c_{2}^{-}}\right)N_{2}\right\}\leq
0.43,
\end{multline*}
and
\begin{multline*}
\max\left\{\frac{\overline{M}_{1}}{\alpha_{1}^{-}},\frac{\overline{M}_{2}}{\alpha_{2}^{-}},\frac{\overline{M}_{3}}{\alpha_{3}^{-}},\left(1+\frac{\alpha_{1}^{+}}{\alpha_{1}^{-}}\right)\overline{M}_{1},
\left(1+\frac{\alpha_{2}^{+}}{\alpha_{2}^{-}}\right)\overline{M}_{2},\left(1+\frac{\alpha_{3}^{+}}{\alpha_{3}^{-}}\right)\overline{M}_{3},\right.\\
\left.\frac{\overline{N}_{1}}{c_{1}^{-}},\frac{\overline{N}_{2}}{c_{2}^{-}},
\left(1+\frac{c_{1}^{+}}{c_{1}^{-}}\right)\overline{N}_{1},\left(1+\frac{c_{2}^{+}}{c_{2}^{-}}\right)\overline{N}_{2}\right\}\leq
1.
\end{multline*}
So, condition $(H_{3})$ holds. Therefore, using Theorem \ref{th0}, \ref{th1} and Theorem \ref{th2}, we conclude that the HOBAMs (\ref{eq1}) with the coefficients and parameters defined above has one and only one Stepanov-like weighted pseudo-almost periodic solution. Besides, this unique solution is globally exponentially stable.\\
Let $h^{\ast}\left( \cdot\right)=\left(x_{1}^{\ast
}\left(\cdot\right),x_{1}^{\ast
}\left(\cdot\right),x_{3}^{\ast}\left(\cdot\right),y_{1}^{\ast}\left(\cdot\right),y_{2}^{\ast}\left(\cdot\right)\right)^{T}$
be a Stepanov-like weighted pseudo-almost periodic on time-space
scales solution of system (\ref{eq1}). We have clearly for $i=1,2,3$
and $j=1,2$
\begin{multline*}
\alpha_{i}^{-}-\sum\limits_{j=1}^{2}\left[L_{j}\left(D_{ij}^{+}+(D_{ij}^{\tau})^{+}
+(\overline{D}_{ij})^{+}\sigma_{ij}^{+}+(\widetilde{D}_{ij})^{+}\xi_{ij}^{+}
+\sum\limits_{k=1}^{2}T_{ijk}^{+}(L_{k}r+|f_{k}(0)|)\right)\right.\\
+\left.(L_{j}r+|f_{j}(0)|)\sum\limits_{k=1}^{2}T_{ijk}^{+}L_{k}\right]>0,
\end{multline*}
and
\begin{multline*}
c_{j}^{-}-\sum\limits_{i=1}^{3}\left[L_{i}\left(E_{ij}^{+}+(E_{ij}^{\tau})^{+}
+(\overline{E}_{ij})^{+}\sigma_{ij}^{+}+(\widetilde{E}_{ij})^{+}\xi_{ij}^{+}
+\sum\limits_{k=1}^{3}T_{ijk}^{+}(L_{k}r+|f_{k}(0)|)\right)\right.\\
+\left.(L_{i}r+|f_{i}(0)|)\sum\limits_{k=1}^{3}T_{ijk}^{+}L_{k}\right]>0,
\end{multline*}
Then from Theorem \ref{th3} all solutions
$\psi=(\varphi_{1},\varphi_{2},\varphi_{3},\phi_{1},\phi_{2})$ of
(\ref{eq1}) satisfying
\begin{equation*}
x_{i}^{\ast }\left( 0\right) =\varphi _{i}\left( 0\right), \,\
y_{j}^{\ast }\left( 0\right) =\phi _{j}\left( 0\right), 1\leq i\leq
3,\,\ 1\leq j\leq 2,
\end{equation*}
converge to its unique Stepanov-like weighted pseudo-almost periodic
on time-space scales solution $h^{\ast}$.
\begin{remark}
The model studied in \cite{maas} is without the second-order
connection weights of delayed feedback.  Moreover, the results
published in \cite{maas} can not be applicable for our example in
this brief. Hence, the results presented here are more general than
that outcomes published in \cite{maas}, \cite{bam0} and
\cite{leak3}. Consequently, this analysis of dynamics behavior of
the Stepanov weighted pseudo almost automorphic on time-space scales
solutions for HOBAMs model with Stepanov-like weighted pseudo almost
automorphic (SWPAA) coefficients and mixed delays improve the
previous study in \cite{maas}, \cite{bam0} and \cite{leak3}.
\end{remark}
%\begin{remark}
%The combination between Stepanov-like weighted pseudo almost
%periodic on time-space scales and almost automorphic on time-space
%scales leads to the Stepanov-like weighted pseudo almost automorphic
%on time-space scales functions, which is considered in this work.
%\end{remark}
%\begin{remark}
%In \cite{aouiti}, the authors investigated the pseudo almost
%automorphic solutions for recurrent neural networks with
%time-varying coefficient and continuously distributed delays. It is
%a good paper and we can used our method, described in this work, to
%extended the pseudo almost automorphic solution studied in
%\cite{aouiti} to the Stepanov-like Weighted Pseudo-Almost
%Automorphic Solutions on Time Scales solutions.
%\end{remark}
\section{Conclusion and open problem}
In this paper, the Stepanov-like weighted pseudo-almost automorphic
on time-space scales are concerned for HOBAMs with mixed delays and
leakage time-varying delays by using Banach's fixed point theorem,
the theory of calculus on time scales and the Lyapunov-Krasovskii
functional method method. The results obtained in this paper are
completely new and complement the previously known works of (Ref.
\cite{maas,bam0,leak3}). Finally, numerical example was given to
demonstrate the effectiveness of our theory. As an interesting
problem, other almost automorphic on time scales type of solutions
for neural networks were deserved to be studied by applying some
suitable methods, such as the study of Stepanov-like weighted
pseudo-almost automorphic on time-space scales for Cohen Grossberg
BAM neural networks. The corresponding results will appear in the
near future.


\begin{thebibliography}{02}
\bibitem{tomorphic2} Wang, C., Ravi P., A., (2014)
\emph{Weighted piecewise pseudo almost automorphic functions with
applications to abstract impulsive $\nabla$-dynamic equations on
time scales}. Advances in Difference Equations 1, 1-29.
\bibitem{chang} Chang, Y. K., Zheng, S., (2016)
\emph{Weighted pseudo almost automorphic solutions to functional
differential equations with infinite delay}. Electronic Journal of
Differential Equations 286, 1-19.
\bibitem{adnene+ahmed} Arbi, A., Alsaedi, A., Cao, J., (2018)
\emph{Delta-differentiable weighted pseudo-almost automorphicity on
time-space scales for a novel class of high-order competitive neural
networks with WPAA coefficients and mixed delays}. Neural Processing
Letters 47, Issue 1, pp 203-232.
\bibitem{carlos} Lizama, C.,  Mesquita, J. G. (2013)
\emph{Almost automorphic solutions of dynamic equations on time
scales}. Journal of Functional Analysis, 265(10), 2267-2311.
\bibitem{diagana} Diagana, T., Mophou, G. M.,  N'Gu\'er\'ekata, G. M. (2010)
\emph{Existence of weighted pseudo-almost periodic solutions to some
classes of differential equations with-weighted pseudo-almost
periodic coefficients}. Nonlinear Analysis: Theory, Methods and
Applications, 72(1), 430-438.
\bibitem{add1} Li, Y., Yang, L., (2014)
\emph{Almost automorphic solution for neutral type high-order
Hopfield neural networks with delays in leakage terms on time
scales}. Applied Mathematics and Computation 242 (2014) 679-693.
\bibitem{add2} Zhu, H., Zhu, Q., Sun, X., Zhou, H.,
\emph{Existence and exponential stability of pseudo almost
automorphic solutions for Cohen-Grossberg neural networks with mixed
delays}. Advances in Difference Equations  (2016) 2016:120.
\bibitem{adn0} Arbi, A., Aouiti, C.,  Touati, A. (2012)
\emph{Uniform asymptotic stability and global asymptotic stability
for time-delay Hopfield neural networks}. Artificial Intelligence
Applications and Innovations, 483-492.
\bibitem{adn00} Arbi, A., Cao, J.,  Alsaedi, A. (2018).
\emph{Improved synchronization analysis of competitive neural networks with time-varying delays}.
NONLINEAR ANALYSIS-MODELLING AND CONTROL, 23(1), 82-102
\bibitem{adn01}Arbi, A., Aouiti, C., Ch\'erif, F., Touati, A.,  Alimi, A. M. (2015)
\emph{Stability analysis of delayed Hopfield neural networks with
impulses via inequality techniques}. Neurocomputing, 158, 281-294.
\bibitem{adn1} Arbi, A., Aouiti, C., Ch\'erif, F., Touati, A., Alimi, A. M., (2015)
\emph{Stability analysis for delayed high-order type of Hopfield
neural networks with impulses}. Neurocomputing, 165, 312-329.
\bibitem{aaa} Arbi, A., Ch\'erif, F., Aouiti, C.,  Touati, A., (2016)
\emph{Dynamics of New Class of Hopfield Neural Networks with
Time-varying and Distributed Delays}. Acta Mathematica Scientia
36(3), 891-912.
\bibitem{adnene+jinde} Arbi, A., Cao, J., (2017)
\emph{Pseudo-almost periodic solution on time-space scales for a
novel class of competitive neutral-type neural networks with mixed
time-varying delays and leakage delays}. Neural Processing Letters,
46 (2), 719-745.
\bibitem{maas} Arbi, A., (2017)
\emph{Dynamics of BAM neural networks with mixed delays and leakage
time-varying delays in the weighted pseudo almost-periodic on
time-space scales}. Mathematical Methods in the Applied Sciences,
1-29.  (In press). DOI: 10.1002/mma.4661.
\bibitem{kren} Ren, K., Qu, J., (2014)
\emph{Identification of Shaft Centerline Orbit for Wind Power Units
Based on Hopfield Neural Network Improved by Simulated Annealing}.
Mathematical Problems in Engineering, 1-6.
\bibitem{zhang1}
Zhang, Y., Li, S., Guo, H. (2017) \emph{A type of biased
consensus-based distributed neural network for path planning}.
Nonlinear Dynamics 89 (3), 1803–1815.
\bibitem{zhang2}
Jin, L., Zhang, Y., Li, S.,  Zhang, Y. (2016) \emph{Modified ZNN for
time-varying quadratic programming with inherent tolerance to noises
and its application to kinematic redundancy resolution of robot
manipulators}. IEEE Transactions on Industrial Electronics, 63(11),
6978-6988.
\bibitem{zhang3}
Zhang, Y.,  Li, S. (2017) \emph{Time-Scale Expansion-Based
Approximated Optimal Control for Underactuated Systems Using
Projection Neural Networks}. IEEE Transactions on Systems, Man, and
Cybernetics: Systems, 1-11. DOI: 10.1109/TSMC.2017.2703140.
\bibitem{kosko} Kosko, B., (1988)
\emph{Bi-directional associative memories}. IEEE Transactions on
Systems, Man and Cybernetics, 18 (1) 49-60.
\bibitem{bam0} Chen, A., Huang, L., Cao, J. (2003)
\emph{Existence and stability of almost periodic solution for BAM
neural networks with delays}. Applied Mathematics and Computation,
137(1), 177-193.
\bibitem{bam1} Berezansky, L., Braverman, E., Idels, L. (2014)
\emph{New global exponential stability criteria for nonlinear delay
differential systems with applications to BAM neural networks}.
Applied Mathematics and Computation, 243, 899-910.
\bibitem{bam2} Cao, J., (2003)
\emph{Global asymptotic stability of delayed bi-directional
associative memory neural networks}. Applied Mathematics and
Computation 142, 333-339.
\bibitem{bam3} Arik, S. Tavsanoglu, V. (2005)
\emph{Global asymptotic stability analysis of bidirectional
associative memory neural networks with constant time delays}.
Neurocomputing 68 161-176.
\bibitem{bam4} Ozcan, N., Arik, S. (2009)
\emph{A new sufficient condition for global robust stability of
bidirectional associative memory neural networks with multiple time
delays}. Nonlinear Analysis: Real World Applications, 10 (5),
3312-3320.
\bibitem{bam5} Huo, H. F., Li, W. T., Tang, S. (2009)
\emph{Dynamics of high-order BAM neural networks with and without
impulses}. Applied Mathematics and Computation, 215 (6), 2120-2133.
\bibitem{bam6} Cao, J., Liang, J., Lam, J. (2004)
\emph{Exponential stability of high-order bidirectional associative
memory neural networks with time delays}. Physica D: Nonlinear
Phenomena, 199(3), 425-436.
\bibitem{ada1} Sowmiya, C., Raja, R., Cao, J., Rajchakit, G.,  Alsaedi, A. (2017).
\emph{Enhanced robust finite-time passivity for Markovian jumping
discrete-time BAM neural networks with leakage delay}. Advances in
Difference Equations, (1), 318.
\bibitem{ada2} Zhang, H., Ye, R., Cao, J.,  Alsaedi, A. (2017).
\emph{Existence and Globally Asymptotic Stability of Equilibrium
Solution for Fractional-Order Hybrid BAM Neural Networks with
Distributed Delays and Impulses}. Complexity.
\bibitem{ada3}
Bao, H.,  Cao, J. (2012). \emph{Exponential stability for stochastic
BAM networks with discrete and distributed delays}. Applied
Mathematics and Computation, 218(11), 6188-6199.
\bibitem{ada4} Jinde Cao, R. Rakkiyappan, K. Maheswari, A.
Chandrasekar, (2016). \emph{Exponential $H_{\infty}$ filtering
analysis for discrete-time switched neural networks with random
delays using sojourn probabilities}, Science China Technological
Sciences, 59:3(2016) 387-402.
\bibitem{highbam1} Zhang, A., Qiu, J., She, J., (2014)
\emph{Existence and global exponential stability of periodic
solution for high-order discrete-time BAM neural networks}. Neural
Networks 50, 98-109.
\bibitem{highbam} Ren, F., Cao, J., (2007)
\emph{Periodic oscillation of higher-order bidirectional associative
memory neural networks with periodic coefficients and delays}.
Nonlinearity 20 605-629.
\bibitem{leak1} Gopalsamy, K. (2007) \emph{Leakage delays in BAM}. Journal of Mathematical Analysis and Applications, 325(2), 1117-1132.
\bibitem{leak2} Liu, B. (2013) \emph{Global exponential stability for BAM neural networks with time-varying delays in the leakage terms}.
Nonlinear Analysis: Real World Applications, 14(1), 559-566.
\bibitem{leak3} Li, Y., Li, Y. (2013)
\emph{Existence and exponential stability of almost periodic
solution for neutral delay BAM neural networks with time-varying
delays in leakage terms}. Journal of the Franklin Institute, 350(9),
2808-2825.
\bibitem{leak0} Gan, Q., Liang, Y. (2012)
\emph{Synchronization of chaotic neural networks with time delay in
the leakage term and parametric uncertainties based on sampled-data
control}. Journal of the Franklin Institute, 349(6), 1955-1971.
\bibitem{hilger} Hilger S. Ein Mabkettenkalk\"{a}ul mit Anwendung auf Zentrumsmannigfaltigkeiten, Ph.D thesis, Universit\"{a}t W\"{a}urzburg, 1988.
\bibitem{advance} Bohner, M., Peterson, A., (2003)
\emph{Advances in Dynamic Equations on Time Scales}. Birkh\"auser
Boston, Boston.
\bibitem{guseinov} Guseinov, G.S., (2003) \emph{Integration on time scales}. J. Math.
Anal. Appl. 285, 107-127.
\bibitem{diaganna} Diagana, T., \emph{Almost Automorphic Type and Almost Periodic Type Functions in Abstract Spaces}. Spring
International Publishing, Switzerland (2013).
\end{thebibliography}
\end{document}